\documentclass[11pt]{amsart}
\usepackage{amsmath,amssymb,mathrsfs}
\newtheorem*{gaussbonnettheorem}{Gauss-Bonnet Theorem}
\newtheorem{theorem}{Theorem}[section]
\newtheorem{corollary}[theorem]{Corollary}

\newtheorem{proposition}[theorem]{Proposition}
\newtheorem{conjecture}[theorem]{Conjecture}
\theoremstyle{definition}
\newtheorem{definition}[theorem]{Definition}
\newcommand{\scal}{\text{\rm scal}}
\newcommand{\Ric}{\text{\rm Ric}{}}
\newcommand{\tracefreeRic}{\overset{\text{\rm o}}{\text{\rm Ric}}{}}

\begin{document}

\title{Curvature, sphere theorems, and the Ricci flow}
\author{Simon Brendle and Richard Schoen}
\address{Department of Mathematics \\ Stanford University \\ Stanford, CA 94305}
\begin{abstract}
In 1926, Hopf proved that any compact, simply connected Riemannian manifold with constant curvature $1$ is isometric to the standard sphere. Motivated by this result, Hopf posed the question whether a compact,
simply connected manifold with suitably pinched curvature is topologically a sphere.

In the first part of this paper, we provide a background discussion, aimed at non-experts, of Hopf's pinching problem and the Sphere Theorem. In the second part, we sketch the proof of the Differentiable Sphere Theorem, and discuss various related results. These results employ a variety of methods, including geodesic and minimal surface techniques as well as Hamilton's Ricci flow.
\end{abstract}
\maketitle 

\section{Hopf's pinching problem and the Sphere Theorem: a background discussion for non-experts}

Differential geometry is concerned with the study of smooth $n$-manifolds. These are topological spaces which are locally homeomorphic to open subsets of $\mathbb{R}^n$ by local coordinate maps such that all change of coordinate maps are smooth diffeomorphisms. The simplest examples of smooth manifolds are two-dimensional surfaces embedded in $\mathbb{R}^3$: these are subsets of $\mathbb{R}^3$ with the property that a neighborhood of each point can be expressed as the graph of a smooth function over some two-dimensional plane. In addition to their topological structure such surfaces inherit two metric structures from $\mathbb{R}^3$. The most elementary metric structure is the one obtained by restricting the distance function from $\mathbb{R}^3$ to $M$. This distance function is called the chord (or extrinsic) distance, and it depends on the embedding of $M$ in $\mathbb{R}^3$. For our purposes, we are interested in the geodesic distance which is obtained by minimizing arc length among all curves on $M$ joining two given endpoints. For example, when we travel along the surface of the earth from one point to another in the most efficient way the distance we travel is the geodesic distance. The geodesic distance is an intrinsic quantity: in order to compute it, one does not need to leave the surface, and so it does not depend on the embedding of the surface into ambient space. We usually encode the geodesic distance by the first fundamental form (or metric tensor) of $M$. This is simply the restriction of the Euclidean inner product to each tangent plane of $M$. From this inner product, we can compute the lengths of tangent vectors, hence the lengths of smooth curves lying on $M$ and, by mimimizing this length, we can find the geodesic distance function. 

Much of the theory of surfaces in $\mathbb{R}^3$ was developed by Gauss in the early 19th century. In addition to the first fundamental form, the theory involves the second fundamental form. Given a point $p \in M$, the second fundamental form is a quadratic form defined on the tangent space to $M$ at $p$. Loosely speaking, the second fundamental form measures the curvature (in $\mathbb{R}^3$) of geodesics in $M$ passing through the point $p$. The eigenvalues of the second fundamental form are called the principal curvatures. If $M$ is convex, then both principal curvatures have the same sign. If the surface $M$ looks like a saddle locally near $p$, then the principal curvatures will have opposite signs. Finally, if $M$ is a plane in $\mathbb{R}^3$, then both principal curvatures vanish.

The product of the principal curvatures defines a real-valued function $K$ on the surface $M$, which is called the Gaussian curvature of $M$. A fundamental discovery of Gauss, which he called his Theorema Egregium, is that the Gaussian curvature of a surface can be computed in terms of the metric tensor and its first two derivatives. Thus this function is uniquely determined by the intrinsic geometry of $M$. Furthermore, Gauss showed that the function $K$ vanishes if and only if $M$ (equipped with the geodesic distance) is locally isometric to the Euclidean plane.

The Theorema Egregium opened up the study of the intrinsic geometry of surfaces. That is, we consider a two-dimensional manifold $M$ equipped with a metric tensor, but we do not require $M$ to be embedded into $\mathbb{R}^3$. By the Theorema Egregium, there is a well-defined notion of Gaussian curvature for such a surface. An important example of such a surface is the real projective plane $\mathbb{RP}^2$ which may be defined as the set of lines through the origin in $\mathbb{R}^3$. Since each line intersects the unit sphere in a pair of antipodal points we may view $\mathbb{RP}^2$ as the surface obtained by identifying pairs of antipodal points on $S^2$. Since the map which sends a point to its antipode is an isometry of $S^2$, we see that $\mathbb{RP}^2$ inherits a metric from $S^2$, and thus is a surface with constant Gaussian curvature $1$. On the other hand $\mathbb{RP}^2$ cannot be isometrically embedded into $\mathbb{R}^3$, so to study projective geometry we need to expand our intuition beyond surfaces in $\mathbb{R}^3$ to the intrinsic geometry of surfaces. 

Of course, most surfaces have non-constant Gaussian curvature. Nevertheless, for any two-dimensional surface the integral of the Gaussian curvature can be expressed in terms of the Euler characteristic of the surface:

\begin{gaussbonnettheorem}
Let $M$ be a compact two-dimensional surface equipped with a metric tensor. Then 
\[\int_M K \, d\text{\rm vol} = 2\pi \chi(M),\] 
where $K$ denotes the Gaussian curvature of $M$, $d\text{\rm vol}$ denotes the induced area measure on $M$ and $\chi(M)$ denotes the Euler characteristic of $M$.
\end{gaussbonnettheorem}

Recall that every compact, orientable surface is diffeomorphic to either $S^2$ or a connected sum of $k$ copies of the torus $T^2$. The number $k$ (which represents the number of handles) is related to the Euler characteristic of $M$ by $\chi(M) = 2-2k$. Therefore, the total curvature of $M$ uniquely determines the Euler characteristic of $M$ which, in turn, determines the diffeomorphism type of $M$. In particular, if $M$ is a compact, orientable surface with positive Gaussian curvature, then $M$ is diffeomorphic to $S^2$. Similarly, if $M$ is a compact, nonorientable surface of positive curvature, then $M$ is diffeomorphic to the real projective plane $\mathbb{RP}^2$. 

In the middle of the 19th century, Riemann extended Gauss' theory of surfaces to higher dimensions. To fix notation, let $M$ denote a compact manifold of dimension $n$ which is equipped with a metric tensor. The metric tensor assigns to each point $p \in M$ a positive definite inner product $g: T_p M \times T_p M \to \mathbb{R}$, which varies smoothly from point to point. Such a structure is referred to as a Riemannian metric (or simply a metric), and the pair $(M,g)$ is referred to as a Riemannian manifold. 

The curvature of a Riemannian manifold $(M,g)$ is described by its Riemann curvature tensor $R$. However, the Riemann curvature tensor of a higher dimensional manifold $(M,g)$ is a much more complicated object than the Gaussian curvature of a surface. A good way to understand it is to consider the sectional curvatures. These have a natural geometric interpretation as the Gaussian curvature of two-dimensional surfaces in $M$. To explain this, we consider all geodesics emanating from $p$ that are tangent to the plane $\pi$. The union of these geodesics rays defines a two-dimensional surface $\Sigma \subset M$; more formally, the surface $\Sigma$ is defined as $\Sigma = \exp_p(U \cap \pi)$, where $\exp_p: T_p M \to M$ denotes the exponential map and $U \subset T_p M$ denotes a small ball centered at the origin. With this understood, the sectional curvature $K(\pi)$ is defined to be the Gaussian curvature of the two-dimensional surface $\Sigma$ at the point $p$. Further discussion of the sectional curvature, and its geometric significance, can be found e.g. in \cite{Gromov2}.

One of the most basic examples of a Riemannian manifold is the $n$-dimensional sphere $S^n$ with its standard metric arising from its embedding as the unit sphere in $\mathbb{R}^{n+1}$. This manifold has constant sectional curvature $1$; that is, $K(\pi) = 1$ for all two-dimensional planes $\pi$. Conversely, it was shown by H.~Hopf in 1926 that a compact, simply connected Riemannian manifold with constant sectional curvature $1$ is necessarily isometric to the sphere $S^n$, equipped with its standard metric (see \cite{Hopf1}, \cite{Hopf2}). More generally, if $(M,g)$ is a compact Riemannian manifold with constant sectional curvature $1$, then $(M,g)$ is isometric to a quotient $S^n / \Gamma$, where $\Gamma$ is a finite group of isometries acting freely. These quotient manifolds are completely classified (see \cite{Wolf}); they are referred to as spherical space forms. The simplest examples of spherical space forms are the sphere $S^n$ and the real projective space $\mathbb{RP}^n$. When $n$ is even, these are the only examples. By contrast, there is an infinite collection of spherical space forms for each odd integer $n$. The lens spaces in dimension $3$ constitute an important class of examples. These spaces were first studied by H.~Tietze \cite{Tietze} in 1908. To describe the definition, let us fix two relatively prime integers $p,q$, and let $\Gamma = \{\alpha \in \mathbb{C}: \alpha^p = 1\}$ denote the cyclic group of $p$-th roots of unity. The group $\Gamma$ acts freely on the unit sphere $S^3 \subset \mathbb{C}^2$ by $\alpha.(z_1,z_2) = (\alpha z_1,\alpha^q z_2)$; hence, the quotient $S^3 / \Gamma$ is a spherical space form.

Hopf conjectured that a compact, simply connected Riemannian manifold whose sectional curvatures are close to $1$ should be homeomorphic to a sphere. (See Marcel Berger's account in \cite{Berger-survey}, page 545.) This idea is formalized by the notion of curvature pinching, which goes back to H.~Hopf and H.E.~Rauch: 

\begin{definition}
\label{global.pinching}
A Riemannian manifold $(M,g)$ is said to be weakly $\delta$-pinched in the global sense if the sectional curvature of $(M,g)$ satisfies $\delta \leq K \leq 1$. If the strict inequality holds, we say that $(M,g)$ is strictly $\delta$-pinched in the global sense.
\end{definition}

For our purposes, it will be convenient to consider the weaker notion of pointwise pinching. This means that we only compare sectional curvatures corresponding to different two-dimensional planes based at the same point $p \in M$:

\begin{definition}
\label{pointwise.pinching}
We say that $(M,g)$ is weakly $\delta$-pinched in the pointwise sense if $0 \leq \delta \, K(\pi_1) \leq K(\pi_2)$ for all points $p \in M$ and all two-dimensional planes $\pi_1,\pi_2 \subset T_p M$. If the strict inequality holds, we say that $(M,g)$ is strictly $\delta$-pinched in the pointwise sense.
\end{definition}

The Differentiable Sphere Theorem, proved in \cite{Brendle-Schoen1}, asserts that any compact Riemannian manifold $(M,g)$ which is strictly $1/4$-pinched in the pointwise sense admits another Riemannian metric which has constant sectional curvature $1$. In particular, this implies that $M$ is diffeomorphic to a spherical space form. In dimension $2$, the Differentiable Sphere Theorem reduces to the statement that a compact surface of positive Gaussian curvature is diffeomorphic to $S^2$ or $\mathbb{RP}^2$. (In dimension $2$, there is only one sectional curvature at each point; hence, every two-dimensional surface of positive curvature is $1/4$-pinched in the pointwise sense.)

The pinching constant $1/4$ is optimal: in fact, any compact symmetric space of rank $1$ admits a metric whose sectional curvatures lie in the interval $[1,4]$. The list of these spaces includes the following examples: 
\begin{itemize}
\item The complex projective space $\mathbb{CP}^m$ (dimension $2m \geq 4$).
\item The quaternionic projective space $\mathbb{HP}^m$ (dimension $4m \geq 8$).
\item The projective plane over the octonions (dimension $16$). 
\end{itemize} 
The manifold $\mathbb{CP}^m$ is defined as the set of complex lines through the origin in $\mathbb{C}^{m+1}$. We may think of the space $\mathbb{CP}^m$ as the quotient of the unit sphere $S^{2m+1} \subset \mathbb{C}^{m+1}$ by the natural $S^1$-action $\alpha.(z_1,\hdots,z_{m+1}) = (\alpha z_1,\hdots,\alpha z_{m+1})$. Since the standard metric on $S^{2m+1}$ is invariant under this $S^1$-action, it induces a natural Riemannian metric on the quotient space $\mathbb{CP}^m = S^{2m+1} / S^1$. 

The manifold $\mathbb{HP}^m$ can similarly be defined as the space of left quaternionic lines through the origin in $\mathbb{H}^{m+1}$. (Note that we need to distinguish left and right quaternionic lines since $\mathbb{H}$ is not commutative.) The space $\mathbb{HP}^m$ can alternatively be characterized as a quotient of $S^{4m+3}$ by a suitable group action. To explain this, we observe that the group of unit quaternions acts freely on the unit sphere $S^{4m+3} \subset \mathbb{H}^{m+1}$ by $\alpha.(z_1,\hdots,z_{m+1}) = (\alpha z_1,\hdots,\alpha z_{m+1})$. Since left multiplication by a unit quaternion preserves the standard metric on $S^{4m+3}$, the quotient space $\mathbb{HP}^m$ inherits a natural Riemannian metric. 

The metrics we have described on $\mathbb{CP}^m$ and $\mathbb{HP}^m$ have sectional curvatures varying between $1$ and $4$ (see \cite{Besse}). For a detailed discussion of the octonions and the octonion projective plane, we refer to the survey article \cite{Baez}.

In addition to the compact symmetric spaces of rank $1$, there are a few other constructions of positively curved metrics. First, the compact homogeneous manifolds of positive sectional curvature have been classified by Berger \cite{Berger4}, Aloff and Wallach \cite{Aloff-Wallach}, and B\'erard-Bergery \cite{Berard-Bergery}. Second, there are biquotient constructions by Eschenburg \cite{Eschenburg1} and Bazaikin \cite{Bazaikin}. The combination of these constructions give additional examples of positively curved metrics in dimensions $6$, $7$, $12$, $13$, and $24$. In dimensions $7$ and $13$, they produce infinitely many topologically distinct examples. 

The strategy used in the proof of the Differentiable Sphere Theorem is to deform the given $1/4$-pinched metric to a metric of constant curvature $1$. This is achieved using Hamilton's Ricci flow. The Ricci flow is a geometric evolution equation of parabolic type; it should be viewed as a nonlinear heat equation for Riemannian metrics. A major obstacle is that, even though the initial metric has $1/4$-pinched curvature, this condition may not be maintained under the evolution. The proof involves finding a suitable curvature condition which is preserved by the Ricci flow and is implied by $1/4$-pinching. This curvature condition is closely related to the notion of positive isotropic curvature, which arises naturally in the context of minimal surface theory! This will be discussed in Section \ref{preserved.curvature.conditions} below.

As pointed out above, the Gauss-Bonnet theorem implies that any compact surface with positive Gaussian curvature is diffeomorphic to $S^2$ or $\mathbb{RP}^2$. More generally, the same argument implies that if $M$ is a compact surface with nonnegative Gaussian curvature, then either $M$ is diffeomorphic to $S^2$ or $\mathbb{RP}^2$, or else the Gaussian curvature of $M$ vanishes identically (in which case $M$ is isometric to a flat torus or Klein bottle). There is a similar borderline statement in the Differentiable Sphere Theorem, which asserts that any compact Riemannian manifold $(M,g)$ which is weakly $1/4$-pinched in the pointwise sense is either diffeomorphic to a spherical space form or isometric to a locally symmetric space. 

In particular, this implies a rigidity theorem. To explain this, let $M$ be a compact symmetric space of rank $1$. As explained above, $M$ admits a canonical metric, whose sectional curvatures lie in the interval $[1,4]$. By the sphere theorem, $M$ does not admit a metric which is strictly $1/4$-pinched in the pointwise sense. Furthermore, any Riemannian metric on $M$ which is weakly $1/4$-pinched in the pointwise sense is isometric to the standard metric up to scaling.

This paper is organized as follows. In Section \ref{topological.sphere.theorem}, we discuss the Topological Sphere Theorem, and the problem of exotic spheres. In Section \ref{variational.techniques}, we discuss some extensions of the Topological Sphere Theorem. This includes the Diameter Sphere Theorem and the Sphere Theorem of Micallef and Moore. These results rely on the variational theory for geodesics and minimal two-spheres, respectively. In Section \ref{ricci.flow}, we describe the Ricci flow method and its applications to Riemannian geometry. It is interesting to compare the different techniques, as each has its own strengths.

In Section \ref{preserved.curvature.conditions}, we sketch some of the key ingredients in the proof of the Differentiable Sphere Theorem. The borderline case in the Differentiable Sphere Theorem is discussed in Section \ref{rigidity.results}. Finally, in Section \ref{other.developments}, we describe various applications to other problems.

\section{The history of the Sphere Theorem and the problem of Exotic Spheres}
\label{topological.sphere.theorem}

To fix notation, let $M$ denote a compact manifold of dimension $n$, and let $g$ be a Riemannian metric on $M$. The curvature of $(M,g)$ is described by the Riemann curvature tensor $R$. This gives, for each point $p \in M$, a multilinear form $R: T_p M \times T_p M \times T_p M \times T_p M \to \mathbb{R}$. Moreover, the Riemann curvature tensor satisfies the symmetries 
\begin{equation} 
\label{alg.1}
R(X,Y,Z,W) = -R(Y,X,Z,W) = R(Z,W,X,Y) 
\end{equation}
and the first Bianchi identity 
\begin{equation}
\label{alg.2}
R(X,Y,Z,W) + R(Y,Z,X,W) + R(Z,X,Y,W) = 0 
\end{equation}
for all tangent vectors $X,Y,Z,W \in T_p M$. By contracting the Riemann curvature tensor, we obtain the Ricci and scalar curvature of $(M,g)$: 
\[\Ric(X,Y) = \sum_{k=1}^n R(X,e_k,Y,e_k)\] 
and 
\[\scal = \sum_{k=1}^n \Ric(e_k,e_k).\] 
Here, $X,Y$ are arbitrary vectors in the tangent space $T_p M$, and $\{e_1,\hdots,e_n\}$ is an orthonormal basis of $T_p M$.

The sectional curvature which was introduced informally in the previous section can now be described precisely: given any point $p \in M$ and any two-dimensional plane $\pi \subset T_p M$, the sectional curvature of $\pi$ is defined by 
\[K(\pi) = \frac{R(X,Y,X,Y)}{|X|^2 \, |Y|^2 - \langle X,Y \rangle^2},\] 
where $\{X,Y\}$ is a basis of $\pi$. Note that this definition is independent of the choice of the basis $\{X,Y\}$.

Hopf's pinching problem was first taken up by H.E.~Rauch after he visited Hopf in Z\"urich during the late 1940s (\cite{Berger-survey}, page 545). In a seminal paper \cite{Rauch}, Rauch showed that a compact, simply connected Riemannian manifold which is strictly $\delta$-pinched in the global sense is homeomorphic to $S^n$ $(\delta \approx 0.75)$. Furthermore, Rauch posed the question of what the optimal pinching constant $\delta$ should be. This question was settled around 1960 by the celebrated Topological Sphere Theorem of M.~Berger and W.~Klingenberg: 

\begin{theorem}[M.~Berger \cite{Berger2}, W.~Klingenberg \cite{Klingenberg}]
\label{top.sphere.thm}
Let $(M,g)$ be a compact, simply connected Riemannian manifold which is strictly $1/4$-pinched in the global sense. Then $M$ is homeomorphic to $S^n$.
\end{theorem}

The classical proof of the Topological Sphere Theorem relies on comparison geometry techniques (see e.g. \cite{Cheeger-Ebin}, Chapter 6). There is an alternative proof due to M.~Gromov, which can be found in \cite{Eschenburg2} (see also \cite{Andrews}).

There are several ways in which one might hope to improve Theorem \ref{top.sphere.thm}. A natural question to ask is whether the global pinching condition in Theorem \ref{top.sphere.thm} can be replaced by a pointwise one. Furthermore, one would like to extend the classification in Theorem \ref{top.sphere.thm} to include manifolds that are not necessarily simply connected. By applying Theorem \ref{top.sphere.thm} to the universal cover, one can conclude that any compact Riemannian manifold which is strictly $1/4$-pinched in the global sense is homeomorphic to a quotient of a sphere by a finite group, but this leaves open the question whether the group acts by standard isometries. We point out that R.~Fintushel and R.J.~Stern \cite{Fintushel-Stern} have constructed an exotic $\mathbb{Z}_2$-action on the standard sphere $S^4$ (see also \cite{Cappell-Shaneson}).

An even more fundamental question is whether a Riemannian manifold satisfying the assumptions of Theorem \ref{top.sphere.thm} is diffeomorphic, instead of just homeomorphic, to $S^n$. This is a highly non-trivial matter, as the smooth structure on $S^n$ is not unique in general. In other words, there exist examples of so-called exotic spheres which are homeomorphic, but not diffeomorphic, to $S^n$. Hence, we may rephrase the problem as follows:

\begin{conjecture}
\label{differentiable.pinching.problem}
An exotic sphere cannot admit a metric with $1/4$-pinched sectional curvature.
\end{conjecture}

The first examples of exotic spheres were constructed in a famous paper by J.~Milnor \cite{Milnor} in 1957. M.~Kervaire and J.~Milnor \cite{Kervaire-Milnor} later proved that there exist exactly $28$ different smooth structures on $S^7$. It was shown by E.~Brieskorn that the exotic $7$-spheres have a natural interpretation in terms of certain affine varieties (cf. \cite{Brieskorn1}, \cite{Brieskorn2}, \cite{Hirzebruch}). To describe this result, let $\Sigma_k$ denote the intersection of the affine variety 
\[\{(z_1,z_2,z_3,z_4,z_5) \in \mathbb{C}^5: z_1^2 + z_2^2 + z_3^2 + z_4^3 + z_5^{6k-1} = 0\}\] 
with the unit sphere in $\mathbb{C}^5$. Brieskorn proved that, for each $k \in \{1,\hdots,28\}$, $\Sigma_k$ is a smooth manifold which is homeomorphic to $S^7$. Moreover, the manifolds $\Sigma_k$, $k \in \{1,\hdots,28\}$, realize all the smooth structures on $S^7$.

In 1974, D.~Gromoll and W.~Meyer \cite{Gromoll-Meyer} described an example of an exotic seven-sphere that admits a metric of nonnegative sectional curvature. It was shown by F.~Wilhelm \cite{Wilhelm} that the Gromoll-Meyer sphere admits a metric which has strictly positive sectional curvature outside a set of measure zero (see also \cite{Eschenburg-Kerin}). P.~Petersen and F.~Wilhelm have recently proposed a construction of a metric of strictly positive sectional curvature on the Gromoll-Meyer sphere, which is currently in the process of verification. 

For each $n \geq 5$, the collection of all smooth structures on $S^n$ has the structure of a finite group $\Theta_n$, called the Kervaire-Milnor group. If $n \equiv 1,2 \mod 8$, there is a natural invariant $\alpha: \Theta_n \to \mathbb{Z}_2$. This invariant is described in more detail in \cite{Joachim-Wraith}. In particular, it is known that half of all smooth structures on $S^n$ have non-zero $\alpha$-invariant. Using the Atiyah-Singer index theorem, N.~Hitchin \cite{Hitchin} showed that an exotic sphere with non-zero $\alpha$-invariant cannot admit a metric of positive scalar curvature.

\begin{theorem}[N.~Hitchin \cite{Hitchin}]
Let $n$ be a positive integer such that either $n \equiv 1 \mod 8$ or $n \equiv 2 \mod 8$. Then half of all smooth structures on $S^n$ do not admit a metric of positive scalar curvature.
\end{theorem}

Conversely, it follows from a theorem of S.~Stolz \cite{Stolz} that every exotic sphere with vanishing $\alpha$-invariant does admit a metric of positive scalar curvature. 

Conjecture \ref{differentiable.pinching.problem} is known as the Differentiable Pinching Problem. This problem has been studied by a large number of authors since the 1960s, and various partial results have been obtained. Gromoll \cite{Gromoll} and Calabi showed that a simply connected Riemannian manifold which is $\delta(n)$-pinched in the global sense is diffeomorphic to $S^n$. The pinching constant $\delta(n)$ depends only on the dimension, and converges to $1$ as $n \to \infty$. In 1971, M.~Sugimoto, K.~Shiohama, and H.~Karcher \cite{Sugimoto-Shiohama} proved an analogous theorem with a pinching constant $\delta$ independent of $n$ $(\delta = 0.87)$. The pinching constant was subsequently improved by E.~Ruh \cite{Ruh1} $(\delta = 0.80)$ and by K.~Grove, H.~Karcher, and E.~Ruh \cite{Grove-Karcher-Ruh2} $(\delta = 0.76)$.

In 1975, H.~Im~Hof and E.~Ruh proved the following theorem, which extends earlier work of Grove, Karcher, and Ruh \cite{Grove-Karcher-Ruh1}, \cite{Grove-Karcher-Ruh2}:

\begin{theorem}[H.~Im~Hof, E.~Ruh \cite{Im-Hof-Ruh}]
There exists a decreasing sequence of real numbers $\delta(n)$ with $\lim_{n \to \infty} \delta(n) = 0.68$ such that the following statement holds: if $M$ is a compact Riemannian manifold of dimension $n$ which is $\delta(n)$-pinched in the global sense, then $M$ is diffeomorphic to a spherical space form.
\end{theorem}

E.~Ruh \cite{Ruh2} has obtained a differentiable version of the sphere theorem under a pointwise pinching condition, albeit with a pinching constant converging to $1$ as $n \to \infty$ (see also Theorem \ref{Huiskens.theorem} below).

In 2007, the authors proved the Differentiable Sphere Theorem with the optimal pinching constant $(\delta = 1/4)$, thereby confirming Conjecture \ref{differentiable.pinching.problem}. 

\begin{theorem}[S.~Brendle, R.~Schoen \cite{Brendle-Schoen1}]
\label{diffeo.sphere.theorem}
Let $(M,g)$ be a compact Riemannian manifold which is strictly $1/4$-pinched in the pointwise sense. Then $M$ is diffeomorphic to a spherical space form. In particular, no exotic sphere admits a metric with strictly $1/4$-pinched sectional curvature.
\end{theorem}

Note that Theorem \ref{diffeo.sphere.theorem} only requires a pointwise pinching condition. (In fact, a much weaker curvature condition suffices; see Theorem \ref{convergence.3} below.)

Finally, there is an analogous pinching problem for manifolds with negative sectional curvature. In other words, if $M$ is a compact Riemannian manifold whose sectional curvatures are sufficiently close to $-1$, can we conclude that $M$ is isometric to a quotient of hyperbolic space by standard isometries? This question was resolved in 1987 by M.~Gromov and W.~Thurston \cite{Gromov-Thurston} (see also \cite{Farrell-Jones}, \cite{Pansu}):

\begin{theorem}[M.~Gromov, W.~Thurston \cite{Gromov-Thurston}] 
Given any integer $n \geq 4$ and any $\delta \in (0,1)$, there exists a compact Riemannian manifold $(M,g)$ of dimension $n$ such that the sectional curvatures of $(M,g)$ lie in the interval $[-1,-\delta]$, but $M$ does not admit a metric of constant negative sectional curvature.
\end{theorem}

\section{Generalizations of the Topological Sphere Theorem}

\label{variational.techniques}

In this section, we describe various results concerning positively curved manifolds. 

\begin{theorem}[M.~Gromov \cite{Gromov1}]
Let $(M,g)$ be a compact Riemannian manifold of dimension $n$ with nonnegative sectional curvature. Then the sum of the Betti numbers of $M$ is bounded by some constant $C(n)$, which depends only on the dimension.
\end{theorem}

K.~Grove and K.~Shiohama \cite{Grove-Shiohama} have obtained an interesting generalization of the Topological Sphere Theorem. In this result, the upper bound for the sectional curvature is replaced by a lower bound for the diameter. We shall sketch an argument due to M.~Berger, which relies on an estimate for the Morse index of geodesics. For a detailed exposition, see \cite{Cheeger-Ebin}, Theorem 6.13, or \cite{Brendle-book}, Theorem 1.15. 

\begin{proposition}
\label{index.of.geodesics}
Let $(M,g)$ be a compact Riemannian manifold of dimension $n$ with sectional curvature $K \geq 1$. Suppose that $p$ and $q$ are two points in $M$ such that $d(p,q) = \text{\rm diam}(M,g) > \frac{\pi}{2}$. Moreover, suppose that $\gamma: [0,1] \to M$ is a geodesic satisfying $\gamma(0) = \gamma(1) = p$. Then $\gamma$ has Morse index at least $n-1$.
\end{proposition}

\textit{Sketch of the proof of Proposition \ref{index.of.geodesics}.} 
Using Toponogov's theorem, one can show that the geodesic $\gamma$ must have length $L(\gamma) > \pi$. Let $I$ denote the index form associated with the second variation of arclength. Then 
\begin{align*} 
I(V,V) &= \int_0^1 \big ( \big | D_{\frac{d}{ds}} V(s) \big |^2 - R(\gamma'(s),V(s),\gamma'(s),V(s)) \big ) \, ds \\ 
&\leq \int_0^1 \big ( \big | D_{\frac{d}{ds}} V(s) \big |^2 - L(\gamma)^2 \, |V(s)|^2 \big ) \, ds 
\end{align*} 
for each vector field $V$ along $\gamma$. 

We next consider the space $\mathscr{H}$ of all vector fields of the form $\sin(\pi s) \, X(s)$, where $X(s)$ is a parallel vector field along $\gamma$. Then $D_{\frac{d}{ds}} D_{\frac{d}{ds}} V(s) = -\pi^2 \, V(s)$ for each vector field $V \in \mathscr{H}$. This implies 
\[I(V,V) \leq (\pi^2 - L(\gamma)^2) \int_0^1 |V(s)|^2 \, ds\] 
for all $V \in \mathscr{H}$. Since $L(\gamma) > \pi$, we conclude that the restriction of $I$ to the vector space $\mathscr{H}$ is negative definite. Since $\dim \mathscr{H} = n-1$, it follows that $\gamma$ has Morse index at least $n-1$. \qed \\

\begin{theorem}[K.~Grove, K.~Shiohama \cite{Grove-Shiohama}] 
\label{diameter.sphere.theorem}
Let $(M,g)$ be a compact Riemannian manifold with sectional curvature $K \geq 1$. If the diameter of $(M,g)$ is greater than $\pi/2$, then $M$ is homeomorphic to $S^n$.
\end{theorem}

\textit{Sketch of the proof of Theorem \ref{diameter.sphere.theorem}.} 
For simplicity, we only consider the case $n \geq 4$. (For $n = 3$, the assertion is a consequence of Theorem \ref{Hamilton.dim.3}.) Let $p$ and $q$ be two points in $M$ of maximal distance, and let $\Omega$ denote the space of all smooth paths $\gamma: [0,1] \to M$ such that $\gamma(0) = \gamma(1) = p$. If $\pi_k(M) \neq 0$ for some $k \in \{1,\hdots,n-1\}$, then there exists a geodesic $\gamma \in \Omega$ with Morse index at most $k-1$. On the other hand, $\gamma$ has Morse index at least $n-1$ by Proposition \ref{index.of.geodesics}. This is a contradiction.

Therefore, we have $\pi_k(M) = 0$ for all $k \in \{1,\hdots,n-1\}$. This implies that $M$ is a homotopy sphere. By work of Freedman \cite{Freedman} and Smale \cite{Smale}, $M$ is homeomorphic to $S^n$. \qed \\

Theorem \ref{diameter.sphere.theorem} requires lower bounds for the sectional curvature and the diameter of $(M,g)$. It is natural to ask whether these assumptions can be replaced by lower bounds for the Ricci curvature and volume of $(M,g)$. An important result in this direction was established by G.~Perelman \cite{Perelman1}. This result was subsequently improved by J.~Cheeger and T.~Colding (cf. \cite{Cheeger-Colding}, Theorem A.1.10; see also \cite{Otsu-Shiohama-Yamaguchi}, Theorem A): 

\begin{theorem}[J.~Cheeger, T.~Colding \cite{Cheeger-Colding}]
For each integer $n \geq 2$, there exists a real number $\delta(n) \in (0,1)$ with the following property: if $(M,g)$ is a compact Riemannian manifold of dimension $n$ with $\Ric \geq (n-1) \, g$ and $\text{\rm vol}(M,g) \geq (1-\delta(n)) \, \text{\rm vol}(S^n(1))$, then $M$ is diffeomorphic to $S^n$.
\end{theorem}

We next describe the Sphere Theorem of Micallef and Moore \cite{Micallef-Moore}. This result improves the Topological Sphere Theorem by weakening the curvature assumptions. To that end, Micallef and Moore introduced a novel curvature condition, which they called positive isotropic curvature: 

\begin{definition}
Let $(M,g)$ be a Riemannian manifold of dimension $n \geq 4$. We say that $(M,g)$ has nonnegative isotropic curvature if 
\begin{align*} 
&R(e_1,e_3,e_1,e_3) + R(e_1,e_4,e_1,e_4) \\ 
&+ R(e_2,e_3,e_2,e_3) + R(e_2,e_4,e_2,e_4) \\ 
&- 2 \, R(e_1,e_2,e_3,e_4) \geq 0 
\end{align*} 
for all points $p \in M$ and all orthonormal four-frames $\{e_1,e_2,e_3,e_4\} \subset T_p M$. Moreover, if the strict inequality holds, we say that $(M,g)$ has positive isotropic curvature.
\end{definition}

For each point $p \in M$, we denote by $T_p^{\mathbb{C}} M = TM \otimes_{\mathbb{R}} \mathbb{C}$ the complexified tangent space to $M$ at $p$. The Riemannian metric $g$ extends to a complex bilinear form $g: T_p^{\mathbb{C}} M \times T_p^{\mathbb{C}} M \to \mathbb{C}$. Similarly, the Riemann curvature tensor extends to a complex multilinear form $R: T_p^{\mathbb{C}} M \times T_p^{\mathbb{C}} M \times T_p^{\mathbb{C}} M \times T_p^{\mathbb{C}} M \to \mathbb{C}$. 

\begin{proposition}
The manifold $(M,g)$ has nonnegative isotropic curvature if and only if $R(\zeta,\eta,\bar{\zeta},\bar{\eta}) \geq 0$ for all points $p \in M$ and all vectors $\zeta,\eta \in T_p^{\mathbb{C}} M$ satisfying $g(\zeta,\zeta) = g(\zeta,\eta) = g(\eta,\eta) = 0$.
\end{proposition}

The key idea of Micallef and Moore is to study harmonic two-spheres instead of geodesics. More precisely, for each map $f: S^2 \to M$, the energy of $f$ is defined by 
\[\mathscr{E}(f) = \frac{1}{2} \int_{S^2} \Big ( \Big | \frac{\partial f}{\partial x} \Big |^2 + \Big | \frac{\partial f}{\partial y} \Big |^2 \Big ) \, dx \, dy,\] 
where $(x,y)$ are the coordinates on $S^2$ obtained by stereographic projection. A map $f: S^2 \to M$ is called harmonic if it is a critical point of the functional $\mathscr{E}(f)$. This is equivalent to saying that 
\[D_{\frac{\partial}{\partial x}} \frac{\partial f}{\partial x} + D_{\frac{\partial}{\partial y}} \frac{\partial f}{\partial y} = 0\] 
at each point on $S^2$. In the special case when $(M,g)$ has positive isotropic curvature, Micallef and Moore obtained a lower bound for the Morse index of harmonic two-spheres:

\begin{proposition}[M.~Micallef, J.D.~Moore \cite{Micallef-Moore}]
\label{index.of.harmonic.two.spheres}
Let $(M,g)$ be a compact Riemannian manifold of dimension $n \geq 4$ with positive isotropic curvature, and let $f: S^2 \to M$ be a nonconstant harmonic map. Then $f$ has Morse index at least $[\frac{n-2}{2}]$.
\end{proposition}

\textit{Sketch of the proof of Theorem \ref{index.of.harmonic.two.spheres}.} 
Let $I: \Gamma(f^*(TM)) \times \Gamma(f^*(TM)) \to \mathbb{R}$ denote the index form associated with the second variation of energy. We may extend $I$ to a complex bilinear form $I: \Gamma(f^*(T^{\mathbb{C}} M)) \times \Gamma(f^*(T^{\mathbb{C}} M)) \to \mathbb{C}$. The complexified index form is given by 
\begin{align} 
\label{second.variation.1}
I(W,\overline{W}) 
&= \int_{S^2} \Big [ g \big ( D_{\frac{\partial}{\partial x}} W,D_{\frac{\partial}{\partial x}} \overline{W} \big ) + g \big ( D_{\frac{\partial}{\partial y}} W,D_{\frac{\partial}{\partial y}} \overline{W} \big ) \Big ] \, dx \, dy \notag \\ 
&- \int_{S^2} \Big [ R \Big ( \frac{\partial f}{\partial x},W,\frac{\partial f}{\partial x},\overline{W} \Big ) + R \Big ( \frac{\partial f}{\partial y},W,\frac{\partial f}{\partial y},\overline{W} \Big ) \Big ] \, dx \, dy 
\end{align} 
for all $W \in \Gamma(f^*(T^{\mathbb{C}} M))$. For abbreviation, let $\frac{\partial f}{\partial z} = \frac{1}{2} \, \big ( \frac{\partial f}{\partial x} - i \, \frac{\partial f}{\partial y} \big )$ and $\frac{\partial f}{\partial \bar{z}} = \frac{1}{2} \, \big ( \frac{\partial f}{\partial x} + i \, \frac{\partial f}{\partial y} \big )$. Moreover, given any section $W \in \Gamma(f^*(T^{\mathbb{C}} M))$, we define 
\[D_{\frac{\partial}{\partial \bar{z}}} W = \frac{1}{2} \, \Big ( D_{\frac{\partial}{\partial x}} W + i \, D_{\frac{\partial}{\partial y}} W \Big ).\] 
With this understood, the formula (\ref{second.variation.1}) can be rewritten as 
\begin{align} 
\label{second.variation.2}
I(W,\overline{W}) 
&= 4 \int_{S^2} g \big ( D_{\frac{\partial}{\partial \bar{z}}} W,D_{\frac{\partial}{\partial z}} \overline{W} \big ) \, dx \, dy \\ 
&- 4 \int_{S^2} R \Big ( \frac{\partial f}{\partial z},W,\frac{\partial f}{\partial \bar{z}},\overline{W} \Big ) \, dx \, dy. \notag
\end{align}
Let 
\[\mathscr{H} = \Big \{ W \in \Gamma(f^*(T^{\mathbb{C}} M)): D_{\frac{\partial}{\partial \bar{z}}} W = 0 \Big \}\] 
denote the space of holomorphic sections of the bundle $f^*(T^{\mathbb{C}} M)$. It follows from the Riemann-Roch theorem that $\dim_{\mathbb{C}} \mathscr{H} \geq n$. Moreover, one can show that $g(\frac{\partial f}{\partial z},W) = 0$ for each section $W \in \mathscr{H}$. Finally, $g(\frac{\partial f}{\partial z},\frac{\partial f}{\partial z}) = 0$ since $f$ is conformal.

To complete the proof, one constructs a subspace $\mathscr{H}_0 \subset \mathscr{H}$ with the following properties: 
\begin{itemize}
\item $\dim_{\mathbb{C}} \mathscr{H}_0 \geq [\frac{n-2}{2}]$.
\item For each section $W \in \mathscr{H}_0$, we have $g(W,W) = 0$.
\item If $W \in \mathscr{H}_0$ and $\frac{\partial f}{\partial z} \wedge W$ vanishes identically, then $W$ vanishes identically.
\end{itemize}
We now consider a section $W \in \mathcal{H}_0$. Since $(M,g)$ has positive isotropic curvature, we have 
\[R \Big ( \frac{\partial f}{\partial z},W,\frac{\partial f}{\partial \bar{z}},\overline{W} \Big ) \geq 0\] 
at each point on $S^2$. This implies 
\[I(W,\overline{W}) = -4 \int_{S^2} R \Big ( \frac{\partial f}{\partial z},W,\frac{\partial f}{\partial \bar{z}},\overline{W} \Big ) \, dx \, dy \leq 0.\] 
Moreover, if $I(W,\overline{W}) = 0$, then $\frac{\partial f}{\partial z} \wedge W$ vanishes identically; consequently, $W$ vanishes identically. 

This shows that the restriction of $I$ to $\mathscr{H}_0$ is negative definite. Thus, we conclude that $\text{\rm ind}(f) \geq \dim_{\mathbb{C}} \mathscr{H}_0 \geq [\frac{n-2}{2}]$, as claimed. \\

Combining Proposition \ref{index.of.harmonic.two.spheres} with the variational theory for harmonic maps (see e.g. \cite{Schoen-Yau}, Chapter VII), Micallef and Moore were able to draw the following conclusion:

\begin{theorem}[M.~Micallef, J.D.~Moore \cite{Micallef-Moore}]
\label{Micallef.Moore.theorem}
Let $(M,g)$ be a compact, simply connected Riemannian manifold of dimension $n \geq 4$ with positive isotropic curvature. Then $M$ is homeomorphic to $S^n$.
\end{theorem}

\textit{Sketch of the proof of Theorem \ref{Micallef.Moore.theorem}.} 
The idea is to study the homotopy groups of $M$. If $\pi_k(M) \neq 0$ for some $k \in \{2,\hdots,[\frac{n}{2}]\}$, then the variational theory for harmonic maps implies that there exists a nonconstant harmonic map $f: S^2 \to M$ with Morse index at most $k - 2$. On the other hand, any nonconstant harmonic map from $S^2$ into $M$ has Morse index at least $[\frac{n-2}{2}]$ by Proposition \ref{index.of.harmonic.two.spheres}. This is a contradiction.

Therefore, we have $\pi_k(M) = 0$ for $k = 2, \hdots, [\frac{n}{2}]$. Since $M$ is assumed to be simply connected, it follows that $H_k(M,\mathbb{Z}) = 0$ for $k = 1, \hdots, [\frac{n}{2}]$. Using Poincar\'e duality, it follows that $H_k(M,\mathbb{Z}) = 0$ for $k = 1, \hdots, n-1$. This shows that $M$ is a homotopy sphere. Hence, it follows from results of Freedman \cite{Freedman} and Smale \cite{Smale} that $M$ is homeomorphic to $S^n$. \qed \\

We note that any manifold $(M,g)$ which is strictly $1/4$-pinched in the pointwise sense has positive isotropic curvature. Hence, Theorem \ref{Micallef.Moore.theorem} generalizes the Topological Sphere Theorem of Berger and Klingenberg.

The topology of non-simply connected manifolds with positive isotropic curvature is not fully understood. It has been conjectured that the fundamental group of a compact manifold $M$ with positive isotropic curvature is virtually free in the sense that it contains a free subgroup of finite index (see e.g. \cite{Fraser}, \cite{Gromov3}). The following result provides some evidence in favor of this conjecture:

\begin{theorem}[A.~Fraser \cite{Fraser}] 
\label{Fraser.thm}
Let $M$ be a compact Riemannian manifold of dimension $n \geq 4$ with positive isotropic curvature. Then the fundamental group of $M$ does not contain a subgroup isomorphic to $\mathbb{Z} \oplus \mathbb{Z}$.
\end{theorem}

The proof of Theorem \ref{Fraser.thm} relies on a delicate analysis of stable minimal tori. This result was proved in dimension $n \geq 5$ by A.~Fraser \cite{Fraser}. In \cite{Brendle-Schoen-survey}, the authors extended Fraser's theorem to the four-dimensional case. The topology of manifolds with positive isotropic curvature is also studied in \cite{Gadgil-Seshadri}.

To conclude this section, we mention some results concerning the topology of four-manifolds with positive sectional curvature (see also \cite{Cao2}). 

\begin{theorem}[W.~Seaman \cite{Seaman}; M.~Ville \cite{Ville}]
\label{seaman.ville.theorem}
Let $(M,g)$ be a compact, orientable Riemannian manifold of dimension $4$ which is $\delta$-pinched in the global sense $(\delta \approx 0.188)$. Then $(M,g)$ is homeomorphic to $S^4$ or $\mathbb{CP}^2$.
\end{theorem}

The proof of Theorem \ref{seaman.ville.theorem} relies on the Bochner formula for harmonic two-forms. This formula had been studied intensively by J.P.~Bourguignon \cite{Bourguignon} and others. The curvature term in the Bochner formula is closely related to the notion of nonnegative isotropic curvature; see e.g. \cite{Brendle-Schoen-survey} or \cite{Micallef-Wang}. In particular, if $(M,g)$ is a compact four-manifold with positive isotropic curvature, then $(M,g)$ does not admit a non-zero harmonic two-form. We note that D.~Meyer \cite{Meyer} has obtained a vanishing theorem for harmonic $k$-forms on manifolds with positive curvature operator. A survey on the Bochner technique can be found in \cite{Berard}.

Finally, let us mention a beautiful theorem of W.Y.~Hsiang and B.~Kleiner concerning positively curved four-manifolds with symmetry.

\begin{theorem}[W.Y.~Hsiang, B.~Kleiner \cite{Hsiang-Kleiner}]
Let $(M,g)$ be a compact, orientable Riemannian manifold of dimension $4$ with positive sectional curvature. If $(M,g)$ admits a non-trivial Killing vector field, then $(M,g)$ is homeomorphic to $S^4$ or $\mathbb{CP}^2$.
\end{theorem}

\section{Deforming Riemannian metrics by the Ricci flow}

\label{ricci.flow}

In this section, we describe the Ricci flow approach. This technique was introduced in seminal work of R.~Hamilton in the 1980s (see e.g. \cite{Hamilton1}, \cite{Hamilton2}). The fundamental idea is to start with a given a Riemannian manifold $(M,g_0)$, and evolve the metric by the evolution equation
\[\frac{\partial}{\partial t} g(t) = -2 \, \Ric_{g(t)}, \qquad g(0) = g_0.\] 
Here, $\Ric_{g(t)}$ denotes the Ricci tensor of the time-dependent metric $g(t)$.

The Ricci flow should be viewed as a heat equation for Riemannian metrics. Hamilton \cite{Hamilton1} proved that, for any choice of the initial metric $g_0$, the Ricci flow has a solution on some maximal time interval $[0,T)$, where $T > 0$ (see also \cite{DeTurck}). Furthermore, if $T < \infty$, then $\liminf_{t \to T} \sup_M |R_{g(t)}| = \infty$. In particular, if the Ricci flow develops a finite-time singularity, then the Riemann curvature tensor of $(M,g(t))$ must be unbounded. This result was later improved by N.~\v Se\v sum \cite{Sesum} who showed that $\limsup_{t \to T} \sup_M |\Ric_{g(t)}| = \infty$ if $T < \infty$.

As an example, suppose that $g_0$ is the standard metric on $S^n$ with constant sectional curvature $1$. In this case, the metrics $g(t) = (1-2(n-1)t) \, g_0$ form a solution to the Ricci flow. This solution is defined for all $t \in [0,\frac{1}{2(n-1)})$, and collapses to a point as $t \to \frac{1}{2(n-1)}$.

In dimension $3$, Hamilton showed that the Ricci flow deforms any initial metric with positive Ricci curvature to a constant curvature metric:

\begin{theorem}[R.~Hamilton \cite{Hamilton1}]
\label{Hamilton.dim.3}
Let $(M,g_0)$ be a compact three-manifold with positive Ricci curvature. Moreover, let $g(t)$, $t \in [0,T)$, denote the unique maximal solution to the Ricci flow with initial metric $g_0$. Then the rescaled metrics $\frac{1}{4(T-t)} \, g(t)$ converge to a metric of constant sectional curvature $1$ as $t \to T$.
\end{theorem} 

The proof of Theorem \ref{Hamilton.dim.3} relies on pointwise curvature estimates. These are established using a suitable version of the maximum principle for tensors.

Theorem \ref{Hamilton.dim.3} has important topological implications. It implies that any compact three-manifold with positive Ricci curvature is diffeomorphic to a spherical space form. Using the classification of spherical space forms in \cite{Wolf}, Hamilton was able to give a complete classification of all compact three-manifolds that admit metrics of positive Ricci curvature.

Hamilton's convergence theorem in dimension $3$ has inspired a large body of work over the last 25 years. In particular, two lines of research have been pursued:

First, one would like to study the global behavior of the Ricci flow in dimension $3$ for general initial metrics (i.e. without the assumption of positive Ricci curvature). This line of research was pioneered by Hamilton, who developed many crucial technical tools (see e.g. \cite{Hamilton4}, \cite{Hamilton5}, \cite{Hamilton-survey}). It culminated in Perelman's proof of the Poincar\'e and Geometrization conjectures (cf. \cite{Perelman2}, \cite{Perelman3}, \cite{Perelman4}). A non-technical survey can be found in \cite{Besson1} or \cite{Leeb}. 

Another natural problem is to extend the convergence theory for the Ricci flow to dimensions greater than $3$. In this case, one assumes that the initial metric satisfies a suitable curvature condition. The goal is to show that the evolved metrics converge to a metric of constant sectional curvature up to rescaling. One of the first results in this direction was established by Hamilton \cite{Hamilton2} in 1986.

\begin{theorem}[R.~Hamilton \cite{Hamilton2}]
\label{Hamilton.dim.4}
Let $(M,g_0)$ be a compact Riemannian manifold of dimension $4$. Assume that $g_0$ has positive curvature operator; that is, $\sum_{i,j,k,l} R_{ijkl} \, \varphi^{ij} \, \varphi^{kl} > 0$ for each point $p \in M$ and every non-zero two-form $\varphi \in \wedge^2 T_p M$. Moreover, let $g(t)$, $t \in [0,T)$, denote the unique maximal solution to the Ricci flow with initial metric $g_0$. Then the rescaled metrics $\frac{1}{6(T-t)} \, g(t)$ converge to a metric of constant sectional curvature $1$ as $t \to T$.
\end{theorem} 

Again, Theorem \ref{Hamilton.dim.4} has a topological corollary: it implies that any compact four-manifold which admits a metric of positive curvature operator is diffeomorphic to $S^4$ or $\mathbb{RP}^4$.

H.~Chen \cite{Chen} proved that the conclusion of Theorem \ref{Hamilton.dim.4} holds under a slightly weaker curvature assumption. A Riemannian manifold $M$ is said to have two-positive curvature operator if $\sum_{i,j,k,l} R_{ijkl} \, (\varphi^{ij} \, \varphi^{kl} + \psi^{ij} \, \psi^{kl}) > 0$ for all points $p \in M$ and all two-forms $\varphi,\psi \in \wedge^2 T_p M$ satisfying $|\varphi|^2 = |\psi|^2 = 1$ and $\langle \varphi,\psi \rangle = 0$. 

\begin{theorem}[H.~Chen \cite{Chen}]
\label{Chen.dim.4}
Let $(M,g_0)$ be a compact Riemannian manifold of dimension $4$ with two-positive curvature operator. Moreover, let $g(t)$, $t \in [0,T)$, denote the unique maximal solution to the Ricci flow with initial metric $g_0$. Then the rescaled metrics $\frac{1}{6(T-t)} \, g(t)$ converge to a metric of constant sectional curvature $1$ as $t \to T$.
\end{theorem} 

Furthermore, Chen \cite{Chen} proved that every four-manifold which is strictly $1/4$-pinched in the pointwise sense has two-positive curvature operator. This is a special feature of the four-dimensional case, which fails in dimension $n \geq 5$. As a consequence, Chen was able to show that every compact four-manifold which is strictly $1/4$-pinched in the pointwise sense is diffeomorphic to $S^4$ or $\mathbb{RP}^4$. B.~Andrews and H.~Nguyen \cite{Andrews-Nguyen} have recently obtained an alternative proof of this result.

We note that C.~Margerin \cite{Margerin2} proved a sharp convergence result for the Ricci flow in dimension $4$. Combining this theorem with techniques from conformal geometry, A.~Chang, M.~Gursky, and P.~Yang proved a beautiful conformally invariant sphere theorem in dimension $4$:

\begin{theorem}[A.~Chang, M.~Gursky, P.~Yang \cite{Chang-Gursky-Yang}]
\label{conformally.invariant.sphere.theorem}
Let $(M,g)$ be a compact four-manifold with positive Yamabe constant. Suppose that $(M,g)$ satisfies the integral pinching condition 
\[\int_M |W|^2 \, d\text{\rm vol} < 16\pi^2 \, \chi(M),\]
where $|W|^2 = \sum_{i,j,k,l} W_{ijkl} \, W^{ijkl}$ denotes the square of the norm of the Weyl tensor of $(M,g)$. Then $M$ is either diffeomorphic to $S^4$ or $\mathbb{RP}^4$.
\end{theorem}

\textit{Sketch of the proof of Theorem \ref{conformally.invariant.sphere.theorem}.} 
Let 
\[\tracefreeRic_g = \Ric_g - \frac{1}{4} \, \scal_g \: g\] 
denote the trace-free part of the Ricci tensor of $(M,g)$. It follows from the Gauss-Bonnet formula that 
\[\int_M \Big ( \frac{1}{6} \, \scal_g^2 - 2 \, |\tracefreeRic_g|^2 + |W_g|^2 \Big ) \, d\text{\rm vol}_g = 32\pi^2 \, \chi(M)\] 
(cf. \cite{Chang-Gursky-Yang}, equation (0.4)). This implies 
\begin{align*} 
&\int_M \Big ( \frac{1}{6} \, \scal_g^2 - 2 \, |\tracefreeRic_g|^2 - |W_g|^2 \Big ) \, d\text{\rm vol}_g \\ 
&= 32\pi^2 \, \chi(M) - 2 \int_M |W_g|^2 \, d\text{\rm vol}_g > 0. 
\end{align*}
The key step in the proof is to construct a conformal metric $\tilde{g} = e^{2w} \, g$ which has positive scalar curvature and satisfies the pointwise inequality 
\[\frac{1}{6} \, \scal_{\tilde{g}}^2 - 2 \, |\tracefreeRic_{\tilde{g}}|^2 - |W_{\tilde{g}}|^2 > 0.\] 
Having constructed a metric $\tilde{g}$ with these properties, a theorem of C.~Margerin \cite{Margerin2} implies that the Ricci flow evolves the metric $\tilde{g}$ to a constant curvature metric. This shows that $M$ is diffeomorphic to either $S^4$ or $\mathbb{RP}^4$. \qed \\

If $g$ is a Riemannian metric on $\mathbb{CP}^2$ which is conformal to the Fubini-Study metric, then 
\[\int_{\mathbb{CP}^2} |W|^2 \, d\text{\rm vol} = 48\pi^2.\] 
Since $\chi(\mathbb{CP}^2) = 3$, this shows that the curvature condition in Theorem \ref{conformally.invariant.sphere.theorem} is sharp.

We next describe some convergence theorems for the Ricci flow that hold in all dimensions. The first result of this type was proved by G.~Huisken in 1985. The following is a corollary of Huisken's theorem:

\begin{theorem}[G.~Huisken \cite{Huisken}]
\label{Huiskens.theorem}
For each integer $n \geq 4$, there exists a real number $\delta(n) \in (0,1)$ with the following property: let $(M,g_0)$ be a compact Riemannian manifold of dimension $n$ which is strictly $\delta(n)$-pinched in the pointwise sense, and let $g(t)$, $t \in [0,T)$, denote the unique maximal solution to the Ricci flow with initial metric $g_0$. Then the rescaled metrics $\frac{1}{2(n-1)(T-t)} \, g(t)$ converge to a metric of constant sectional curvature $1$ as $t \to T$.
\end{theorem} 

We note that C.~Margerin \cite{Margerin1} and S.~Nishikawa \cite{Nishikawa} have also obtained convergence results for the Ricci flow in arbitrary dimension. 

The following result extends Theorem \ref{Chen.dim.4} to higher dimensions:

\begin{theorem}[C.~B\"ohm, B.~Wilking \cite{Bohm-Wilking}]
Let $(M,g_0)$ be a compact Riemannian manifold with two-positive curvature operator, and let $g(t)$, $t \in [0,T)$, denote the unique maximal solution to the Ricci flow with initial metric $g_0$. Then the rescaled metrics $\frac{1}{2(n-1)(T-t)} \, g(t)$ converge to a metric of constant sectional curvature $1$ as $t \to T$.
\end{theorem}

\section{Curvature conditions that are preserved by the Ricci flow}

\label{preserved.curvature.conditions}

All known convergence theorems for the Ricci flow share some common features. In particular, they all exploit the fact that a certain curvature condition is preserved by the Ricci flow. In this section, we describe some general tools for verifying that a given curvature condition is preserved by the Ricci flow. These tools are based on the maximum principle, and were developed by Hamilton \cite{Hamilton1}, \cite{Hamilton2}. 

Let $g(t)$, $t \in [0,T)$, be a solution to the Ricci flow on a manifold $M$. Moreover, let $E$ denote the pull-back of the tangent bundle $TM$ under the map 
\[M \times (0,T) \to M, \quad (p,t) \mapsto p.\] 
Clearly, $E$ is a vector bundle over $M \times (0,T)$, and the fiber of $E$ over the point $(p,t) \in M \times (0,T)$ is given by the tangent space $T_p M$. The sections of the vector bundle $E$ can be viewed as vector fields on $M$ that vary in time. Given any section $X$ of $E$, we define the covariant time derivative of $X$ by 
\[D_{\frac{\partial}{\partial t}} X = \frac{\partial}{\partial t} X - \sum_{k=1}^n \Ric(X,e_k) \, e_k,\] 
where $\{e_1,\hdots,e_n\}$ is a local orthonormal frame with respect to the metric $g(t)$. The covariant time derivative $D_{\frac{\partial}{\partial t}}$ is metric compatible in the sense that 
\begin{align*} 
\frac{\partial}{\partial t} (g(X,Y)) 
&= g \big ( \frac{\partial}{\partial t} X,Y \big ) + g \big ( X,\frac{\partial}{\partial t} Y \big ) - 2 \, \Ric(X,Y) \\ 
&= g \big ( D_{\frac{\partial}{\partial t}} X,Y \big ) + g \big ( X,D_{\frac{\partial}{\partial t}} Y \big ) 
\end{align*}
for all sections $X,Y$ of the bundle $E$.

The Riemann curvature tensor of $g(t)$ can now be viewed as a section of the bundle $E^* \otimes E^* \otimes E^* \otimes E^*$. Furthermore, the covariant time derivative on $E$ induces a covariant time derivative on the bundle $E^* \otimes E^* \otimes E^* \otimes E^*$. With this understood, the evolution equation of the Riemann curvature tensor can be written in the form  
\[D_{\frac{\partial}{\partial t}} R = \Delta R + Q(R),\] 
where $D_{\frac{\partial}{\partial t}}$ denotes the covariant time derivative and $Q(R)$ is a quadratic expression in the Riemannian curvature tensor $R$: 
\begin{align}
\label{definition.of.Q}
Q(R)(X,Y,Z,W) &= \sum_{p,q=1}^n R(X,Y,e_p,e_q) \, R(Z,W,e_p,e_q) \notag \\ 
&+ 2 \sum_{p,q=1}^n R(X,e_p,Z,e_q) \, R(Y,e_p,W,e_q) \\ 
&- 2 \sum_{p,q=1}^n R(X,e_p,W,e_q) \, R(Y,e_p,Z,e_q). \notag
\end{align} 
This evolution equation was first derived by Hamilton \cite{Hamilton2}; see also \cite{Brendle-book}, Section 2.3. The use of the covariant time derivative was suggested by Uhlenbeck.

We next describe Hamilton's maximum principle for the Ricci flow. To fix notation, let $\mathscr{C}_B(\mathbb{R}^n)$ denote the space of all algebraic curvature operators on $\mathbb{R}^n$. In other words, $\mathscr{C}_B(\mathbb{R}^n)$ consists of all multilinear forms $R: \mathbb{R}^n \times \mathbb{R}^n \times \mathbb{R}^n \times \mathbb{R}^n \to \mathbb{R}$ satisfying the relations 
\[R(X,Y,Z,W) = -R(Y,X,Z,W) = R(Z,W,X,Y)\] 
and 
\[R(X,Y,Z,W) + R(Y,Z,X,W) + R(Z,X,Y,W) = 0\] 
for all vectors $X,Y,Z,W \in \mathbb{R}^n$. Moreover, let $F$ be a subset of $\mathscr{C}_B(\mathbb{R}^n)$ which is invariant under the natural action of $O(n)$. Since $F$ is $O(n)$-invariant, it makes sense to say that the curvature tensor of a Riemannian manifold $(M,g)$ lies in the set $F$. To explain this, we fix a point $p \in M$. After identifying the tangent space $T_p M$ with $\mathbb{R}^n$, we may view the curvature tensor of $(M,g)$ at $p$ as an element of $\mathscr{C}_B(\mathbb{R}^n)$. Of course, the identification of $T_p M$ with $\mathbb{R}^n$ is not canonical, but this does not cause problems since $F$ is $O(n)$-invariant.

\begin{theorem}[R.~Hamilton \cite{Hamilton2}]
\label{Hamilton.maximum.principle}
Let $F \subset \mathscr{C}_B(\mathbb{R}^n)$ be a closed, convex set which is invariant under the natural action of $O(n)$. Moreover, we assume that $F$ is invariant under the ODE $\frac{d}{dt} R = Q(R)$. Finally, let $g(t)$, $t \in [0,T)$, be a solution to the Ricci flow on a compact manifold $M$ with the property that the curvature tensor of $(M,g(0))$ lies in $F$ for all points $p \in M$. Then the curvature tensor of $(M,g(t))$ lies in $F$ for all points $p \in M$ and all $t \in [0,T)$.
\end{theorem}

In the remainder of this section, we discuss some important examples of curvature conditions that are preserved by the Ricci flow. Hamilton \cite{Hamilton2} proved that nonnegative curvature operator is preserved in all dimensions. Furthermore, Hamilton showed that nonnegative Ricci curvature is preserved by the Ricci flow in dimension $3$, and nonnegative isotropic curvature is preserved in dimension $4$ (see \cite{Hamilton1}, \cite{Hamilton7}). It turns out that nonnegative Ricci curvature is not preserved by the Ricci flow in dimension $n \geq 4$ (see \cite{Maximo}). By contrast, nonnegative isotropic curvature is preserved by the Ricci flow in all dimensions:

\begin{theorem}[S.~Brendle, R.~Schoen \cite{Brendle-Schoen1}; H.~Nguyen \cite{Nguyen}]
\label{pic.is.preserved}
Let $M$ be a compact manifold of dimension $n \geq 4$, and let $g(t)$, $t \in [0,T)$, be a solution to the Ricci flow on $M$. If $(M,g(0))$ has nonnegative isotropic curvature, then $(M,g(t))$ has nonnegative isotropic curvature for all $t \in [0,T)$.
\end{theorem}

The proof of Theorem \ref{pic.is.preserved} requires two ingredients: the first is Hamilton's maximum principle for the Ricci flow (cf. Theorem \ref{Hamilton.maximum.principle}); the second one is an algebraic inequality for curvature tensors with nonnegative isotropic curvature. We give a sketch of the proof here. A complete proof can be found in \cite{Brendle-book}, Sections 7.2 and 7.3.

\begin{proposition}
\label{pic.is.preserved.2}
Let $R$ be an algebraic curvature tensor on $\mathbb{R}^n$ with nonnegative isotropic curvature. Moreover, suppose that 
\[R_{1313} + R_{1414} + R_{2323} + R_{2424} - 2 \, R_{1234} = 0\] 
for some orthonormal four-frame $\{e_1,e_2,e_3,e_4\}$. Then 
\[Q(R)_{1313} + Q(R)_{1414} + Q(R)_{2323} + Q(R)_{2424} - 2 \, Q(R)_{1234} \geq 0.\]
\end{proposition}

\textit{Sketch of the proof of Proposition \ref{pic.is.preserved.2}.} 
Using the definition of $Q(R)$, we compute 
\begin{align*} 
&Q(R)_{1313} + Q(R)_{1414} + Q(R)_{2323} + Q(R)_{2424} - 2 \, Q(R)_{1234} \\ 
&= \sum_{p,q=1}^n (R_{13pq} - R_{24pq})^2 + \sum_{p,q=1}^n (R_{14pq} + R_{23pq})^2 \\ 
&+ 2 \sum_{p,q=1}^n (R_{1p1q} + R_{2p2q}) \, (R_{3p3q} + R_{4p4q}) - 2 \sum_{p,q=1}^n R_{12pq} \, R_{34pq} \\ 
&- 2 \sum_{p,q=1}^n (R_{1p3q} + R_{2p4q}) \, (R_{3p1q} + R_{4p2q}) \\ 
&- 2 \sum_{p,q=1}^n (R_{1p4q} - R_{2p3q}) \, (R_{4p1q} - R_{3p2q}). 
\end{align*} 
Hence, it suffices to prove that 
\begin{align}
\label{key.inequality}
&\sum_{p,q=1}^n (R_{1p1q} + R_{2p2q}) \, (R_{3p3q} + R_{4p4q}) - \sum_{p,q=1}^n R_{12pq} \, R_{34pq} \notag \\ 
&\geq \sum_{p,q=1}^n (R_{1p3q} + R_{2p4q}) \, (R_{3p1q} + R_{4p2q}) \\ 
&+ \sum_{p,q=1}^n (R_{1p4q} - R_{2p3q}) \, (R_{4p1q} - R_{3p2q}). \notag
\end{align} 
In order to prove (\ref{key.inequality}), we view the isotropic curvature as a real-valued function defined on the space of all orthonormal four-frames. By assumption, this function attains its global minimum at the point $\{e_1,e_2,e_3,e_4\}$. Hence, the first variation at the point $\{e_1,e_2,e_3,e_4\}$ is zero, and the second variation is nonnegative. In order to take advantage of this information, we consider three different types of variations:

Step 1: We first consider the orthonormal four-frame $\{e_1,\cos(s) \, e_2 - \sin(s) \, e_3,\sin(s) \, e_2 + \cos(s) \, e_3,e_4\}$. Since the first variation of the isotropic curvature is zero, we have 
\[R_{1213} + R_{1242} + R_{3413} + R_{3442} = 0.\] 
An analogous argument gives 
\[R_{1214} + R_{1223} + R_{3414} + R_{3423} = 0.\] 
Using these identities, one can show that 
\begin{align}
\label{step.1}
&\sum_{p,q=1}^4 (R_{1p1q} + R_{2p2q}) \, (R_{3p3q} + R_{4p4q}) - \sum_{p,q=1}^4 R_{12pq} \, R_{34pq} \notag \\ 
&= \sum_{p,q=1}^4 (R_{1p3q} + R_{2p4q}) \, (R_{3p1q} + R_{4p2q}) \\ 
&+ \sum_{p,q=1}^4 (R_{1p4q} - R_{2p3q}) \, (R_{4p1q} - R_{3p2q}). \notag
\end{align} 

Step 2: We next consider the four-frame $\{\cos(s) \, e_1 + \sin(s) \, e_q,e_2,e_3,e_4\}$, where $q \in \{5,\hdots,n\}$. Since the first variation of the isotropic curvature is equal to zero, it follows that 
\[R_{133q} + R_{144q} + R_{432q} = 0.\] 
Using this and other analogous identities, we obtain 
\begin{align}
\label{step.2}
&\sum_{p=1}^4 (R_{1p1q} + R_{2p2q}) \, (R_{3p3q} + R_{4p4q}) - \sum_{p=1}^4 R_{12pq} \, R_{34pq} \notag \\ 
&= \sum_{p=1}^4 (R_{1p3q} + R_{2p4q}) \, (R_{3p1q} + R_{4p2q}) \\ 
&+ \sum_{p=1}^4 (R_{1p4q} - R_{2p3q}) \, (R_{4p1q} - R_{3p2q}) \notag
\end{align} 
for $q = 5, \hdots, n$.

Step 3: To describe the third type of variation, we consider four vectors $w_1,w_2,w_3,w_4 \in \text{\rm span}\{e_5,\hdots,e_n\}$. For each $i \in \{1,2,3,4\}$, we denote by $v_i(s)$ the unique solution of the linear ODE 
\[v_i'(s) = \sum_{j=1}^4 (\langle v_i(s),e_j \rangle \, w_j - \langle v_i(s),w_j \rangle \, e_j)\] 
with initial condition $v_i(0) = e_i$. Then $v_i'(0) = w_i$. Moreover, it is easy to see that the vectors $\{v_1(s),v_2(s),v_3(s),v_4(s)\}$ are orthonormal for all $s \in \mathbb{R}$. Since the second variation of the isotropic curvature is nonnegative, we conclude that 
\begin{align} 
\label{second.variation}
0 
&\leq R(w_1,e_3,w_1,e_3) + R(w_1,e_4,w_1,e_4) \notag \\ 
&+ R(w_2,e_3,w_2,e_3) + R(w_2,e_4,w_2,e_4) \notag \\ 
&+ R(e_1,w_3,e_1,w_3) + R(e_2,w_3,e_2,w_3) \notag \\ 
&+ R(e_1,w_4,e_1,w_4) + R(e_2,w_4,e_2,w_4) \notag \\ 
&- 2 \, \big [ R(e_3,w_1,e_1,w_3) + R(e_4,w_1,e_2,w_3) \big ] \\ 
&- 2 \, \big [ R(e_4,w_1,e_1,w_4) - R(e_3,w_1,e_2,w_4) \big ] \notag \\ 
&+ 2 \, \big [ R(e_4,w_2,e_1,w_3) - R(e_3,w_2,e_2,w_3) \big ] \notag \\ 
&- 2 \, \big [ R(e_3,w_2,e_1,w_4) + R(e_4,w_2,e_2,w_4) \big ] \notag \\ 
&- 2 \, R(w_1,w_2,e_3,e_4) - 2 \, R(e_1,e_2,w_3,w_4) \notag
\end{align} 
for all vectors $w_1,w_2,w_3,w_4 \in \text{\rm span}\{e_5,\hdots,e_n\}$. We next define linear transformations $A,B,C,D,E,F: \text{\rm span}\{e_5,\hdots,e_n\} \to \text{\rm span}\{e_5,\hdots,e_n\}$ by 
\[\begin{array}{l@{\qquad}l} 
\langle Ae_p,e_q \rangle = R_{1p1q} + R_{2p2q}, & \langle Be_p,e_q \rangle = R_{3p3q} +  R_{4p4q}, \\ 
\langle Ce_p,e_q \rangle = R_{3p1q} + R_{4p2q}, & \langle De_p,e_q \rangle = R_{4p1q} - R_{3p2q}, \\ 
\langle Ee_p,e_q \rangle = R_{12pq}, & \langle Fe_p,e_q \rangle = R_{34pq} \end{array}\] 
for $p,q \in \{5,\hdots,n\}$. The inequality (\ref{second.variation}) implies that the symmetric operator 
\[\begin{bmatrix} B & F & -C^* & -D^* \\ -F & B & D^* & -C^* \\ -C & D & A & E \\ -D & -C & -E & A \end{bmatrix}\] 
is positive semi-definite. From this, we deduce that 
\[\text{\rm tr}(AB) + \text{\rm tr}(EF) \geq \text{\rm tr}(C^2) + \text{\rm tr}(D^2),\] 
hence 
\begin{align}
\label{step.3}
&\sum_{p,q=5}^n (R_{1p1q} + R_{2p2q}) \, (R_{3p3q} + R_{4p4q}) - \sum_{p,q=5}^n R_{12pq} \, R_{34pq} \notag \\ 
&\geq \sum_{p,q=5}^n (R_{1p3q} + R_{2p4q}) \, (R_{3p1q} + R_{4p2q}) \\ 
&+ \sum_{p,q=5}^n (R_{1p4q} - R_{2p3q}) \, (R_{4p1q} - R_{3p2q}). \notag
\end{align} 

Combining (\ref{step.1}), (\ref{step.2}), and (\ref{step.3}), the inequality (\ref{key.inequality}) follows. \qed \\

We next describe various curvature conditions that are related to nonnegative isotropic curvature, and are also preserved by the Ricci flow. The following is an immediate consequence of Theorem \ref{pic.is.preserved}:

\begin{corollary}[S.~Brendle, R.~Schoen \cite{Brendle-Schoen1}]
\label{product.with.flat.factors}
Let $M$ be a compact manifold of dimension $n \geq 4$, and let $g(t)$, $t \in [0,T)$, be a solution to the Ricci flow on $M$. Then: 
\begin{itemize}
\item If $(M,g(0)) \times \mathbb{R}$ has nonnegative isotropic curvature, the product $(M,g(t)) \times \mathbb{R}$ has nonnegative isotropic curvature for all $t \in [0,T)$. 
\item If $(M,g(0)) \times \mathbb{R}^2$ has nonnegative isotropic curvature, the product $(M,g(t)) \times \mathbb{R}^2$ has nonnegative isotropic curvature for all $t \in [0,T)$.
\end{itemize}
\end{corollary}

Another result in this direction was proved by the first author in \cite{Brendle1} (see also \cite{Brendle-book}, Section 7.6). In the following, $S^2(1)$ denotes a two-dimensional sphere of constant curvature $1$.

\begin{theorem}[S.~Brendle \cite{Brendle1}] 
\label{product.with.S2}
Let $M$ be a compact manifold of dimension $n \geq 4$, and let $g(t)$, $t \in [0,T)$, be a solution to the Ricci flow on $M$. If $(M,g(0)) \times S^2(1)$ has nonnegative isotropic curvature, then $(M,g(t)) \times S^2(1)$ has nonnegative isotropic curvature for all $t \in [0,T)$. 
\end{theorem}

Unlike Corollary \ref{product.with.flat.factors}, Theorem \ref{product.with.S2} does not follow directly from Theorem \ref{pic.is.preserved}. This is because the manifolds $(M,g(t)) \times S^2(1)$ do not form a solution to the Ricci flow.

We now discuss the product conditions in more detail. To that end, we assume that $(M,g)$ is a Riemannian manifold of dimension $n$. We first consider the case $n = 3$. In this case, the following holds:
\begin{itemize}
\item The product $(M,g) \times \mathbb{R}$ has nonnegative isotropic curvature if and only if $(M,g)$ has nonnegative Ricci curvature.
\item The product $(M,g) \times \mathbb{R}^2$ has nonnegative isotropic curvature if and only if $(M,g)$ has nonnegative sectional curvature.
\end{itemize}
For $n \geq 4$, the following proposition gives a necessary and sufficient condition for the product $(M,g) \times \mathbb{R}$ to have nonnegative isotropic curvature.

\begin{proposition}
\label{characterization.of.pic.1}
Let $(M,g)$ be a Riemannian manifold of dimension $n \geq 4$. Then the following statements are equivalent: 
\begin{itemize}
\item[(i)] The product $(M,g) \times \mathbb{R}$ has nonnegative isotropic curvature.
\item[(ii)] We have 
\begin{align*} 
&R(e_1,e_3,e_1,e_3) + \lambda^2 \, R(e_1,e_4,e_1,e_4) \\ 
&+ R(e_2,e_3,e_2,e_3) + \lambda^2 \, R(e_2,e_4,e_2,e_4) \\ 
&- 2\lambda \, R(e_1,e_2,e_3,e_4) \geq 0 
\end{align*}
for all points $p \in M$, all orthonormal four-frames $\{e_1,e_2,e_3,e_4\} \subset T_p M$, and all $\lambda \in [0,1]$. 
\item[(iii)] We have $R(\zeta,\eta,\bar{\zeta},\bar{\eta}) \geq 0$ for all points $p \in M$ and all vectors $\zeta,\eta \in T_p^{\mathbb{C}} M$ satisfying $g(\zeta,\zeta) \, g(\eta,\eta) - g(\zeta,\eta)^2 = 0$.
\end{itemize}
\end{proposition}

The proof of Proposition \ref{characterization.of.pic.1} is purely algebraic (for details, see \cite{Brendle-book}, Proposition 7.18). We next consider the condition that $(M,g) \times \mathbb{R}^2$ has nonnegative isotropic curvature (cf. \cite{Brendle-book}, Proposition 7.18). 

\begin{proposition}
\label{characterization.of.pic.2}
Let $(M,g)$ be a Riemannian manifold of dimension $n \geq 4$. Then the following statements are equivalent: 
\begin{itemize}
\item[(i)] The product $(M,g) \times \mathbb{R}^2$ has nonnegative isotropic curvature.
\item[(ii)] We have 
\begin{align*} 
&R(e_1,e_3,e_1,e_3) + \lambda^2 \, R(e_1,e_4,e_1,e_4) \\ 
&+ \mu^2 \, R(e_2,e_3,e_2,e_3) + \lambda^2\mu^2 \, R(e_2,e_4,e_2,e_4) \\ 
&- 2\lambda\mu \, R(e_1,e_2,e_3,e_4) \geq 0 
\end{align*}
for all points $p \in M$, all orthonormal four-frames $\{e_1,e_2,e_3,e_4\} \subset T_p M$, and all $\lambda,\mu \in [0,1]$. 
\item[(iii)] We have $R(\zeta,\eta,\bar{\zeta},\bar{\eta}) \geq 0$ for all points $p \in M$ and all vectors $\zeta,\eta \in T_p^{\mathbb{C}} M$.
\end{itemize}
\end{proposition}

In the special setting of K\"ahler geometry, Bando \cite{Bando} and Mok \cite{Mok} proved that nonnegative holomorphic bisectional curvature is preserved by the K\"ahler-Ricci flow. Cao and Hamilton \cite{Cao-Hamilton} showed that the weaker notion of nonnegative orthogonal bisectional curvature is preserved by the K\"ahler-Ricci flow as well. The latter condition is related to the notion of nonnegative isotropic curvature: in fact, any K\"ahler manifold with nonnegative isotropic curvature necessarily has nonnegative orthogonal bisectional curvature.

Theorem \ref{pic.is.preserved} and Corollary \ref{product.with.flat.factors} provide important examples of preserved curvature conditions. Each of these curvature conditions defines a closed, convex, $O(n)$-invariant cone in $\mathscr{C}_B(\mathbb{R}^n)$, which is preserved by the Hamilton ODE. By adapting a technique of B\"ohm and Wilking \cite{Bohm-Wilking}, it is possible to construct a family of so-called pinching cones, which are all preserved by the Hamilton ODE. Combining these ideas with general results of R.~Hamilton (see \cite{Hamilton2} or \cite{Brendle-book}, Section 5.4), one can draw the following conclusion: 

\begin{theorem}[S.~Brendle, R.~Schoen \cite{Brendle-Schoen1}]
\label{convergence.1}
Let $(M,g_0)$ be a compact Riemannian manifold of dimension $n \geq 4$ with the property that 
\begin{align*} 
&R(e_1,e_3,e_1,e_3) + \lambda^2 \, R(e_1,e_4,e_1,e_4) \\ 
&+ \mu^2 \, R(e_2,e_3,e_2,e_3) + \lambda^2\mu^2 \, R(e_2,e_4,e_2,e_4) \\ 
&- 2\lambda\mu \, R(e_1,e_2,e_3,e_4) > 0 
\end{align*}
for all orthonormal four-frames $\{e_1,e_2,e_3,e_4\}$ and all $\lambda,\mu \in [0,1]$. Let $g(t)$, $t \in [0,T)$, denote the unique maximal solution to the Ricci flow with initial metric $g_0$. Then the rescaled metrics $\frac{1}{2(n-1)(T-t)} \, g(t)$ converge to a metric of constant sectional curvature $1$ as $t \to T$. 
\end{theorem}

It turns out that any Riemannian manifold of dimension $n \geq 4$ which is strictly $1/4$-pinched in the pointwise sense satisfies the assumption of Theorem \ref{convergence.1}. Hence, we can draw the following conclusion:

\begin{corollary}[S.~Brendle, R.~Schoen \cite{Brendle-Schoen1}] 
\label{convergence.2}
Let $(M,g_0)$ be a compact Riemannian manifold of dimension $n \geq 4$ which is strictly $1/4$-pinched in the pointwise sense. Let $g(t)$, $t \in [0,T)$, denote the unique maximal solution to the Ricci flow with initial metric $g_0$. Then the rescaled metrics $\frac{1}{2(n-1)(T-t)} \, g(t)$ converge to a metric of constant sectional curvature $1$ as $t \to T$. 
\end{corollary}

\textit{Sketch of the proof of Corollary \ref{convergence.2}.} 
For each point $p \in M$, we denote by $K_{\text{\rm max}}(p)$ and $K_{\text{\rm min}}(p)$ the maximum and minimum of the sectional curvature of $(M,g_0)$ at the point $p$. Since $(M,g_0)$ is strictly $1/4$-pinched, we have $0 < K_{\text{\rm max}}(p) < 4 \, K_{\text{\rm min}}(p)$ for each point $p \in M$. It follows from an inequality of Berger that the curvature tensor of $(M,g_0)$ satisfies 
\[R(e_1,e_2,e_3,e_4) \leq \frac{2}{3} \, (K_{\text{\rm max}}(p) - K_{\text{\rm min}}(p))\] 
for all orthonormal four-frames $\{e_1,e_2,e_3,e_4\} \subset T_p M$. This implies 
\begin{align*} 
&R(e_1,e_3,e_1,e_3) + \lambda^2 \, R(e_1,e_4,e_1,e_4) \\ 
&+ \mu^2 \, R(e_2,e_3,e_2,e_3) + \lambda^2\mu^2 \, R(e_2,e_4,e_2,e_4) \\ 
&- 2\lambda\mu \, R(e_1,e_2,e_3,e_4) \\ 
&\geq (1 + \lambda^2 + \mu^2 + \lambda^2 \mu^2) \, K_{\text{\rm min}}(p) - \frac{4}{3} \, \lambda\mu \, (K_{\text{\rm max}}(p) - K_{\text{\rm min}}(p)) \\ 
&= \big ( (1-\lambda\mu)^2 + (\lambda-\mu)^2 \big ) \, K_{\text{\rm min}}(p) + \frac{4}{3} \, \lambda\mu \, (4 \, K_{\text{\rm min}}(p) - K_{\text{\rm max}}(p)) \\ 
&> 0
\end{align*}
for all orthonormal four-frames $\{e_1,e_2,e_3,e_4\} \subset T_p M$ and all $\lambda,\mu \in [0,1]$. Hence, the assertion follows from Theorem \ref{convergence.1}. \qed \\

The Differentiable Sphere Theorem (Theorem \ref{diffeo.sphere.theorem} above) is an immediate consequence of Corollary \ref{convergence.2}.

To conclude this section, we state an improved convergence result for the Ricci flow:

\begin{theorem}[S.~Brendle \cite{Brendle1}]
\label{convergence.3}
Let $(M,g_0)$ be a compact Riemannian manifold of dimension $n \geq 4$ with the property that 
\begin{align*} 
&R(e_1,e_3,e_1,e_3) + \lambda^2 \, R(e_1,e_4,e_1,e_4) \\ 
&+ R(e_2,e_3,e_2,e_3) + \lambda^2 \, R(e_2,e_4,e_2,e_4) \\ 
&- 2\lambda \, R(e_1,e_2,e_3,e_4) > 0 
\end{align*}
for all orthonormal four-frames $\{e_1,e_2,e_3,e_4\}$ and all $\lambda \in [0,1]$. Let $g(t)$, $t \in [0,T)$, denote the unique maximal solution to the Ricci flow with initial metric $g_0$. Then the rescaled metrics $\frac{1}{2(n-1)(T-t)} \, g(t)$ converge to a metric of constant sectional curvature $1$ as $t \to T$. 
\end{theorem}

Theorem \ref{convergence.3} extends many known convergence results for the Ricci flow (see also \cite{Besson2}). The main ingredient in the proof is Theorem \ref{product.with.S2}. A detailed argument can be found in \cite{Brendle-book}, Section 8.4.

\section{The borderline case in the Sphere Theorem}

\label{rigidity.results}

In this section, we describe various rigidity results. We will only sketch the main ideas involved in the proofs. For a detailed exposition, we refer to \cite{Brendle-book}, Chapter 9.

The first result in this direction was established by M.~Berger \cite{Berger3} (see also \cite{Cheeger-Ebin}, Theorem 6.6).

\begin{theorem}[M.~Berger \cite{Berger3}]
Let $(M,g)$ be a compact, simply connected Riemannian manifold which is weakly $1/4$-pinched in the global sense. Then $M$ is either homeomorphic to $S^n$ or isometric to a symmetric space.
\end{theorem}

The borderline case in the Diameter Sphere was studied by D.~Gromoll and K.~Grove \cite{Gromoll-Grove} (see also \cite{Cao-Tang}). 

We now describe some rigidity results obtained by means of the Ricci flow. The following result was established by R.~Hamilton \cite{Hamilton2}: 

\begin{theorem}[R.~Hamilton \cite{Hamilton2}]
\label{rigidity.dim.3}
Let $(M,g_0)$ be a compact three-manifold which is locally irreducible and has nonnegative Ricci curvature. Moreover, let $g(t)$, $t \in [0,T)$, denote the unique maximal solution to the Ricci flow with initial metric $g_0$. Then the rescaled metrics $\frac{1}{4(T-t)} \, g(t)$ converge to a metric of constant sectional curvature $1$ as $t \to T$.
\end{theorem} 

\textit{Sketch of the proof of Theorem \ref{rigidity.dim.3}.} 
By assumption, $(M,g_0)$ is locally irreducible. Hence, if we choose $\tau \in (0,T)$ sufficiently small, then $(M,g(\tau))$ is locally irreducible as well. Furthermore, the manifold $(M,g(\tau))$ has nonnegative Ricci curvature. 

We now consider the following subset of the tangent bundle $TM$: 
\[\{v \in TM: \Ric_{g(\tau)}(v,v) = 0\}.\] 
It follows from the strict maximum principle that this set is invariant under parallel transport with respect to the metric $g(\tau)$. Since $(M,g(\tau))$ is locally irreducible, any parallel subbundle of $TM$ has rank $0$ or $3$. From this, we deduce that 
\[\{v \in T_p M: \Ric_{g(\tau)}(v,v) = 0\} = \{0\}\] 
for each point $p \in M$. Consequently, the manifold $(M,g(\tau))$ has positive Ricci curvature. The assertion now follows from Theorem \ref{Hamilton.dim.3}. \qed \\

In dimension $n \geq 4$, we have the following result:  

\begin{theorem}[S.~Brendle, R.~Schoen \cite{Brendle-Schoen2}]
\label{strict.maximum.principle.for.pic}
Let $M$ be a compact manifold of dimension $n \geq 4$, and let $g(t)$, $t \in [0,T]$ be a solution to the Ricci flow on $M$ with nonnegative isotropic curvature. Then, for each $\tau \in (0,T)$, the set of all orthonormal four-frames $\{e_1,e_2,e_3,e_4\}$ satisfying 
\begin{align*} 
&R_{g(\tau)}(e_1,e_3,e_1,e_3) + R_{g(\tau)}(e_1,e_4,e_1,e_4) \\ 
&+ R_{g(\tau)}(e_2,e_3,e_2,e_3) + R_{g(\tau)}(e_2,e_4,e_2,e_4) \\ 
&- 2 \, R_{g(\tau)}(e_1,e_2,e_3,e_4) = 0 
\end{align*} 
is invariant under parallel transport with respect to the metric $g(\tau)$.
\end{theorem}

In particular, if the reduced holonomy group of $(M,g(\tau))$ is $\text{\rm SO}(n)$, then $(M,g(\tau))$ has positive isotropic curvature.

Theorem \ref{strict.maximum.principle.for.pic} is similar in spirit to a result of R.~Hamilton \cite{Hamilton2} concerning solutions to the Ricci flow with nonnegative curvature operator. However, Hamilton's techniques are not applicable in this setting. Instead, the proof of Theorem \ref{strict.maximum.principle.for.pic} relies on a variant of J.M.~Bony's strict maximum principle for degenerate elliptic equations (cf. \cite{Bony}). This technique was first employed in the context of geometric flows in \cite{Brendle-Schoen2}. It has since found applications to other borderline situations involving Ricci flow (see e.g. \cite{Andrews-Nguyen}, \cite{Gu}). 

Theorem \ref{strict.maximum.principle.for.pic} is particularly effective in combination with M.~Berger's classification of holonomy groups (see \cite{Berger1} or \cite{Besse}, Corollary 10.92). For example, the following structure theorem for compact, simply connected manifolds with nonnegative isotropic curvature was established in \cite{Brendle3}:

\begin{theorem}[S.~Brendle \cite{Brendle3}] 
\label{weakly.pic}
Let $(M,g_0)$ be a compact, simply connected Riemannian manifold of dimension $n \geq 4$ which is irreducible and has nonnegative isotropic curvature. Then one of the following statements holds:
\begin{itemize}
\item[(i)] $M$ is homeomorphic to $S^n$.
\item[(ii)] $n = 2m$ and $(M,g_0)$ is a K\"ahler manifold.
\item[(iii)] $(M,g_0)$ is isometric to a symmetric space.
\end{itemize}
\end{theorem}

\textit{Sketch of the proof of Theorem \ref{weakly.pic}.} 
Suppose that $(M,g_0)$ is non-symmetric. Then there exists a real number $\delta > 0$ such that $(M,g(t))$ is irreducible and non-symmetric for all $t \in (0,\delta)$. By Berger's holonomy theorem, there are three possibilities: 

\textit{Case 1:} Suppose that $\text{\rm Hol}(M,g(\tau)) = \text{\rm SO}(n)$ for some $\tau \in (0,\delta)$. In this case, Theorem \ref{strict.maximum.principle.for.pic} implies that $(M,g(\tau))$ has positive isotropic curvature. Consequently, $M$ is homeomorphic to $S^n$ by Theorem \ref{Micallef.Moore.theorem}.

\textit{Case 2:} Suppose that $n = 2m$ and $\text{\rm Hol}(M,g(t)) = \text{\rm U}(m)$ for all $t \in (0,\delta)$. Then $(M,g(t))$ is K\"ahler for all $t \in (0,\delta)$. Consequently, $(M,g_0)$ is also K\"ahler.

\textit{Case 3:} Suppose that $n = 4m \geq 8$ and $\text{\rm Hol}(M,g(\tau)) = \text{\rm Sp}(m) \cdot \text{\rm Sp}(1)$ for some $\tau \in (0,\delta)$. In this case, $(M,g(\tau))$ is quaternionic-K\"ahler. By Corollary 15 in \cite{Brendle3}, $(M,g(\tau))$ is isometric to $\mathbb{HP}^m$ up to scaling. This contradicts the fact that $(M,g(\tau))$ is non-symmetric. \qed \\

M.~Berger \cite{Berger7} has shown that any quaternionic-K\"ahler manifold with positive sectional curvature is isometric to $\mathbb{HP}^m$ up to scaling. More recently, H.~Seshadri proved that any K\"ahler manifold which satisfies the assumptions of Theorem \ref{weakly.pic} is biholomorphic to complex projective space or isometric to a Hermitian symmetric space (see \cite{Seshadri1}, Theorem 1.2). The proof is similar in spirit to Siu and Yau's proof of the Frankel conjecture (cf. \cite{Siu-Yau}). 

We now state another consequence of Theorem \ref{strict.maximum.principle.for.pic}:

\begin{theorem} 
\label{borderline.case.1}
Let $(M,g_0)$ be a compact, locally irreducible Riemannian manifold of dimension $n \geq 4$ with the property that 
\begin{align*} 
&R(e_1,e_3,e_1,e_3) + \lambda^2 \, R(e_1,e_4,e_1,e_4) \\ 
&+ R(e_2,e_3,e_2,e_3) + \lambda^2 \, R(e_2,e_4,e_2,e_4) \\ 
&- 2\lambda \, R(e_1,e_2,e_3,e_4) \geq 0 
\end{align*}
for all orthonormal four-frames $\{e_1,e_2,e_3,e_4\}$ and all $\lambda \in [0,1]$. Moreover, let $g(t)$, $t \in [0,T)$, denote the unique maximal solution to the Ricci flow with initial metric $g_0$. Then one of the following statements holds: 
\begin{itemize}
\item[(i)] The rescaled metrics $\frac{1}{2(n-1)(T-t)} \, g(t)$ converge to a metric of constant sectional curvature $1$ as $t \to T$.
\item[(ii)] $n = 2m$ and the universal cover of $(M,g_0)$ is a K\"ahler manifold.
\item[(iii)] $(M,g_0)$ is locally symmetric.
\end{itemize}
\end{theorem}

We now impose the stronger assumption that $(M,g_0)$ is weakly $1/4$-pinched in the pointwise sense. In this case, $(M,g_0)$ satisfies the curvature assumption in Theorem \ref{borderline.case.1}. Moreover, if $(M,g_0)$ is K\"ahler, then $(M,g_0)$ is isometric to complex projective space up to scaling. Hence, we can draw the following conclusion: 

\begin{corollary}[S.~Brendle, R.~Schoen \cite{Brendle-Schoen2}] 
\label{borderline.case.2}
Let $(M,g_0)$ be a compact Riemannian manifold of dimension $n \geq 4$ which is weakly $1/4$-pinched in the pointwise sense. Moreover, we assume that $(M,g_0)$ is not locally symmetric. Let $g(t)$, $t \in [0,T)$, denote the unique maximal solution to the Ricci flow with initial metric $g_0$. Then the rescaled metrics $\frac{1}{2(n-1)(T-t)} \, g(t)$ converge to a metric of constant sectional curvature $1$ as $t \to T$.
\end{corollary}

In the remainder of this section, we describe some results concerning almost $1/4$-pinched manifolds. The first result in this direction was proved by M.~Berger in 1983:

\begin{theorem}[M.~Berger \cite{Berger8}]
\label{berger.almost.quarter.pinching}
For every even integer $n$, there exists a real number $\delta(n) \in (0,1/4)$ with the following property: if $(M,g_0)$ is a compact, simply connected Riemannian manifold of dimension $n$ which is strictly $\delta(n)$-pinched in the global sense, then $M$ is homeomorphic to $S^n$ or diffeomorphic to a compact symmetric space of rank one.
\end{theorem}

The proof of Theorem \ref{berger.almost.quarter.pinching} is by contradiction, and relies on a compactness argument in the spirit of Gromov. In particular, the value of the pinching constant $\delta(n)$ is not known in general.

U.~Abresch and W.~Meyer \cite{Abresch-Meyer} showed that any compact, simply connected, odd-dimensional Riemannian manifold whose sectional curvatures lie in the interval $(\frac{1}{4(1+10^{-6})^2},1]$ is homeomorphic to a sphere.

Using the classification in Theorem \ref{borderline.case.1} and a Cheeger-Gromov-style compactness argument, P.~Petersen and T.~Tao obtained the following result:

\begin{theorem}[P.~Petersen, T.~Tao \cite{Petersen-Tao}]
\label{almost.quarter.pinching}
For each integer $n \geq 4$, there exists a real number $\delta(n) \in (0,1/4)$ with the following property: if $(M,g_0)$ is a compact, simply connected Riemannian manifold of dimension $n$ which is strictly $\delta(n)$-pinched in the global sense, then $M$ is diffeomorphic to a sphere or a compact symmetric space of rank one.
\end{theorem}

The conclusion of Theorem \ref{almost.quarter.pinching} can be improved slightly when $n$ is odd. In this case, there exists a real number $\delta(n) \in (0,1/4)$ with the property that every compact $n$-dimensional manifold $(M,g_0)$ which is strictly $\delta(n)$-pinched in the global sense is diffeomorphic to a spherical space form.

\section{Other recent developments}

\label{other.developments}

In this final section, we describe some other applications of the techniques described in Sections \ref{ricci.flow} -- \ref{rigidity.results}. 

A much-studied problem in Riemannian geometry is to classify all Einstein manifolds satisfying a suitable curvature condition. This question was first studied by M.~Berger \cite{Berger5}, \cite{Berger6} in the 1960. Berger showed that if $(M,g)$ is a compact Einstein manifold of dimension $n$ which is strictly $\frac{3n}{7n-4}$-pinched in the global sense, then $(M,g)$ has constant sectional curvature (see \cite{Besse}, Section 0.33). In 1974, S.~Tachibana \cite{Tachibana} proved that any compact Einstein manifold with positive curvature operator has constant sectional curvature. Furthermore, Tachibana showed that a compact Einstein manifold with nonnegative curvature operator is locally symmetric. Other results in this direction were obtained M.~Gursky and C.~LeBrun \cite{Gursky-LeBrun} and D.~Yang \cite{Yang}. The following result provides a classification of Einstein metrics with nonnegative isotropic curvature (see also \cite{Micallef-Wang}, where the four-dimensional case is discussed):

\begin{theorem}[S.~Brendle \cite{Brendle3}, \cite{Brendle4}]
\label{generalization.of.Tachibana}
Let $(M,g)$ be a compact Einstein manifold of dimension $n \geq 4$. If $(M,g)$ has positive isotropic curvature, then $(M,g)$ is isometric to a spherical space form. Furthermore, if $(M,g)$ has nonnegative isotropic curvature, then $(M,g)$ is locally symmetric.
\end{theorem}

H.~Seshadri \cite{Seshadri2} has obtained a classification of almost-Einstein manifolds with nonnegative isotropic curvature. The proof combines Theorem \ref{generalization.of.Tachibana} with a Cheeger-Gromov-style compactness argument.

We next describe some results concerning ancient solutions to the Ricci flow. Recall that a solution $(M,g(t))$ to the Ricci flow is said to be ancient if it is defined for all $t \in (-\infty,0)$ (cf. \cite{Hamilton5}). P.~Daskalopoulos, R.~Hamilton, and N.~\v Se\v sum \cite{Daskalopoulos-Hamilton-Sesum} have recently obtained a complete classification of all ancient solutions to the Ricci flow in dimension $2$. V.~Fateev \cite{Fateev} constructed an interesting example of an ancient solution in dimension $3$ (see also \cite{Bakas-Kong-Ni}).

In \cite{Brendle-Huisken-Sinestrari}, it was shown that any ancient solution to the Ricci flow in dimension $n \geq 3$ which satisfies a suitable curvature pinching condition must have constant sectional curvature. In particular, in dimension $3$ the following result holds:

\begin{theorem}[S.~Brendle, G.~Huisken, C.~Sinestrari \cite{Brendle-Huisken-Sinestrari}]
\label{ancient.dim.3}
Let $M$ be a compact three-manifold, and let $g(t)$, $t \in (-\infty,0)$, be an ancient solution to the Ricci flow on $M$. Moreover, suppose that there exists a uniform constant $\rho > 0$ such that 
\[\Ric_{g(t)} \geq \rho \, \scal_{g(t)} \, g(t) \geq 0\] 
for all $t \in (-\infty,0)$. Then the manifold $(M,g(t))$ has constant sectional curvature for each $t \in (-\infty,0)$.
\end{theorem}

The proof of Theorem \ref{ancient.dim.3} relies on a new interior estimate for the Ricci flow in dimension $3$. To describe this estimate, suppose that $g(t)$, $t \in [0,T)$, is a solution to the Ricci flow on a compact three-manifold $M$. Moreover, suppose that there exists a uniform constant $\rho > 0$ such that 
\[\Ric_{g(t)} \geq \rho \, \scal_{g(t)} \, g(t) \geq 0\] 
for each $t \in [0,T)$. Then, for each $t \in (0,T)$, the curvature tensor of $(M,g(t))$ satisfies the pointwise estimate 
\[|\tracefreeRic_{g(t)}|^2 \leq \Big ( \frac{3}{2t} \Big )^\sigma \, \scal_{g(t)}^{2-\sigma},\] 
where $\sigma = \rho^2$.

In dimension $n \geq 4$, we have the following generalization of Theorem \ref{ancient.dim.3}: 

\begin{theorem}[S.~Brendle, G.~Huisken, C.~Sinestrari \cite{Brendle-Huisken-Sinestrari}]
\label{ancient.higher.dim}
Let $M$ be a compact manifold of dimension $n \geq 4$, and let $g(t)$, $t \in (-\infty,0)$, be an ancient solution to the Ricci flow on $M$. Moreover, suppose that there exists a uniform constant $\rho > 0$ with the following property: for each $t \in (-\infty,0)$, the curvature tensor of $(M,g(t))$ satisfies 
\begin{align*} 
&R_{g(t)}(e_1,e_3,e_1,e_3) + \lambda^2 \, R_{g(t)}(e_1,e_4,e_1,e_4) \\ 
&+ R_{g(t)}(e_2,e_3,e_2,e_3) + \lambda^2 \, R_{g(t)}(e_2,e_4,e_2,e_4) \\ 
&- 2\lambda \, R_{g(t)}(e_1,e_2,e_3,e_4) \geq \rho \, \scal_{g(t)} \geq 0 
\end{align*}
for all orthonormal four-frames $\{e_1,e_2,e_3,e_4\}$ and all $\lambda \in [0,1]$. Then the manifold $(M,g(t))$ has constant sectional curvature for each $t \in (-\infty,0)$.
\end{theorem}

We next discuss an important gradient estimate for the scalar curvature along the Ricci flow. R.~Hamilton first established such an inequality in dimension $2$; see \cite{Hamilton3}. Later, Hamilton \cite{Hamilton4} proved a more general inequality, known the differential Harnack inequality, which holds for any solution to the Ricci flow with nonnegative curvature operator (see also \cite{Cao1}). It turns out that a similar inequality holds under the much weaker condition that $M \times \mathbb{R}^2$ has nonnegative isotropic curvature:

\begin{theorem}[S.~Brendle \cite{Brendle2}]
Let $(M,g(t))$, $t \in [0,T)$, be a solution to the Ricci flow which is complete and has bounded curvature. Moreover, suppose that the product $(M,g(t)) \times \mathbb{R}^2$ has nonnegative isotropic curvature for each $t \in (0,T)$. Then 
\[\frac{\partial}{\partial t} \text{\rm scal} + \frac{1}{t} \, \text{\rm scal} + 2 \, \partial_i \text{\rm scal} \, v^i + 2 \, \text{\rm Ric}(v,v) \geq 0\] 
for all points $(p,t) \in M \times (0,T)$ and all vectors $v \in T_p M$. 
\end{theorem}

The differential Harnack inequality for the Ricci flow has interesting implications for complete manifolds with pointwise pinched curvature. These results are motivated in part by a theorem of Hamilton \cite{Hamilton6} which asserts that any complete, strictly convex hypersurface in $\mathbb{R}^n$ with pointwise pinched second fundamental form is compact. The following result can be viewed as an intrinsic analogue of Hamilton's theorem:

\begin{theorem}[B.~Chen, X.~Zhu \cite{Chen-Zhu}]
\label{Chen.Zhu.theorem}
Given any integer $n \geq 4$, there exists a real number $\delta(n) \in (0,1)$ with the following property: if $(M,g_0)$ is a complete Riemannian manifold of dimension $n$ which has bounded curvature and is strictly $\delta(n)$-pinched in the pointwise sense, then $M$ is compact.
\end{theorem}

Another result in this direction was proved in \cite{Ni-Wu}. 

Theorem \ref{Chen.Zhu.theorem} can be generalized as follows: 

\begin{theorem}[S.~Brendle, R.~Schoen \cite{Brendle-Schoen-survey}]
\label{generalization.of.chen.zhu}
Let $(M,g_0)$ be a complete Riemannian manifold of dimension $n \geq 4$ with bounded curvature. We assume that there exists a uniform constant $\rho > 0$ such that 
\begin{align*} 
&R(e_1,e_3,e_1,e_3) + \lambda^2 \, R(e_1,e_4,e_1,e_4) \\ 
&+ \mu^2 \, R(e_2,e_3,e_2,e_3) + \lambda^2\mu^2 \, R(e_2,e_4,e_2,e_4) \\ 
&- 2\lambda\mu \, R(e_1,e_2,e_3,e_4) \geq \rho \, \scal > 0 
\end{align*} 
for all orthonormal four-frames $\{e_1,e_2,e_3,e_4\}$ and all $\lambda,\mu \in [-1,1]$. Then $M$ is compact.
\end{theorem}

\end{document}